\documentclass{elsarticle}
\usepackage{color}

\usepackage[T1,T2A]{fontenc}	
\usepackage[utf8]{inputenc}	
\usepackage[russian,british]{babel}

\usepackage[a4paper, margin=50pt]{geometry}

\usepackage{subcaption}

 \usepackage{amsmath}
 \usepackage{amssymb}
 \usepackage{amsthm}
\usepackage{xparse}
\usepackage{enumitem}
\usepackage{mathtools} 
\usepackage{tikz-cd}
\usepackage{csquotes}

\usepackage[normalem]{ulem}
\usepackage{natbib}
\usepackage{bm}
\usepackage{dsfont}

\usepackage{graphicx}
\graphicspath{{images/}}

\usepackage{tabularx}

\usepackage{hyperref}
\usepackage{xcolor}
\usepackage{tcolorbox}
\hypersetup{%
	colorlinks=true,
	linkcolor=blue!66!black,
	citecolor=red!66!black,
	urlcolor=green!33!black
}

\theoremstyle{definition}
\newtheorem{thm}{Theorem}[section]
\newtheorem{lemma}[thm]{Lemma}

\newtheorem{defn}[thm]{Definition}
\newtheorem{prop}[thm]{Proposition}
\newtheorem*{hypothesis}{Conjecture}

\newtheorem{remark}{Remark}[section]

\newtheorem{construction}[thm]{Construction}
\renewcommand{\l}{\lambda}
\newcommand\G{\Gamma}

\renewcommand{\emptyset}{\varnothing}

\usepackage{scalerel}

\usepackage{tikz}
\usetikzlibrary{hobby}
\newcommand{\eye}{%
	\kern0.1ex%
	\begin{tikzpicture}[use Hobby shortcut]%
		\draw[thick,line cap=round] (0ex,0.5ex) .. (1ex,1ex) .. (2ex,0.5ex);%
		\draw[thick,line cap=round] (0ex,0.5ex) .. (1ex,0ex) .. (2ex,0.5ex);%
		\filldraw[thick,line cap=round] (1ex,0.5ex) circle (0.15ex);%
	\end{tikzpicture}%
	\kern0.1ex%
}

\usepackage{extarrows}

\newcommand\seq{\Longrightarrow}

\renewcommand{\nabla}{\triangledown}

\let\ophi\phi
\let\ovarphi\varphi
\renewcommand{\phi}{\ovarphi}
\renewcommand{\varphi}{\ophi}

\newcommand\N{\mathbb{N}}

\renewcommand\C{\mathbb{C}}

\renewcommand\sl{\mathfrak{sl}}

\newcommand\corank{\operatorname{corank}}

\usepackage{float}

\newcommand{\setp}[1]{\left\{#1\right\}}

\makeatletter

\def\env@sqcases{%
	\let\@ifnextchar\new@ifnextchar
	\left\lbrack
	\def\arraystretch{1.2}%
	\array{@{}l@{\quad}l@{}}%
}
\makeatother

\ExplSyntaxOn
\NewDocumentCommand{\colvec}{O{,} m}{
	\vector_main:nnnn{\\}{#1}{#2}
}
\seq_new:N \l__vector_arg_seq
\cs_new_protected:Npn \vector_main:nnnn #1 #2 #3{
	\seq_set_split:Nnn \l__vector_arg_seq{#2}{#3}
	\begin{pmatrix}
		\seq_use:Nnnn \l__vector_arg_seq{#1}{#1}{#1}
	\end{pmatrix}
}
\ExplSyntaxOff

\setlist[itemize]{leftmargin=2cm}

\usepackage{xstring}
\usepackage{ifthen}
\usepackage{etoolbox}
\usetikzlibrary{
	knots,
	hobby,
	calc,
	shapes.geometric,
	plotmarks,
	arrows.meta,
	tqft,
	fit,
	patterns,
	patterns.meta,
	decorations.markings,
	fadings,
	decorations.pathmorphing,
	decorations.pathreplacing,
	fixedpointarithmetic,
    external
}
\makeatletter
\tikzset{
	>={Latex[round,width=2mm,length=2mm]},
	l/.style={
		line cap=round,
		line join=round
	},
	c/.style={
		circle,
		inner sep=0pt,
		outer sep=0pt,
		minimum size=#1,
	},
	c/.default=1cm,
	show direction/.style={
		postaction=decorate,
		decoration={
			markings,
			mark=at position #1 with {\arrow[opacity=1]{>>}}
		}
	},
	tcenter/.style={
		baseline={([yshift=#1]current bounding box.center)}
	},
	tcenter/.default=-0.5ex,
    nb/.style={
        draw,
        fill=white,
        thick,
        inner sep=0pt,
        minimum size=#1
    },
    nb/.default=10pt,
    nbe/.style={
        draw,
        fill=white,
        thick,
        inner sep=0pt,
        text width=#1,
        align=center,
        minimum width=#1,
        minimum height=11pt
    },
    nbe/.default=22pt
}

\pgfdeclaredecoration{width to zero}{initial}{
	\state{initial}[width=0pt, next state=line, persistent precomputation={%
		\pgfmathdivide{\pgflinewidth}{\pgfdecoratedpathlength}%
		\let\speed=\pgfmathresult%
		\pgfmathmultiply{\speed}{.5pt}%
		\let\decrement=\pgfmathresult%
		\let\curwidth=\pgflinewidth%
 	}]{}
	\state{line}[width=.5pt, persistent postcomputation={%
		\pgfmathsubtract{\curwidth}{\decrement}%
		\let\curwidth=\pgfmathresult%
	}]{%
		\pgfsetlinewidth{\curwidth}%
		\pgfsetarrows{-}%
		\pgfpathmoveto{\pgfpointorigin}%
		\pgfpathlineto{\pgfqpoint{.75pt}{0pt}}%
		\pgfusepath{stroke}%
	}
	\state{final}{%
		\pgfsetlinewidth{\pgflinewidth}%
		\pgfpathmoveto{\pgfpointorigin}%
		\pgfusepath{stroke}% 
	}
}
\tikzset{%
	vanish/.style={%
		fixed point arithmetic,
		decoration={%
			snake,
			segment length=#1,
			amplitude=0.6mm
		},
		decorate,
		postaction={%
			draw,
			decorate,
			decoration={width to zero}
		}
	},
	vanish/.default=6mm,
	red pseudo/.style={%
			line width=0.75pt,
			decoration={%
				ticks,
				segment length=1.5pt,
				amplitude=0.8pt
			},
			decorate
	}
}

\makeatother

\usepackage{fp}
\def\chordR{1}
\def\chordRad{1cm}
\def\chordCircle{(0,0) circle (\chordRad)}

\def\chordPointDSet#1#2#3{\coordinate (#1) at ($(0,0)!\chordR!90-360*#3/#2:(1,0)$)}
\DeclareRobustCommand{\chord}[3]{%
    \tikzset{external/export next=false}%
	\begin{tikzpicture}[scale=#1,tcenter]
		\draw[thick, black, l] \chordCircle;
		\foreach \i/\j in {#3} {
			\chordPointDSet{a}{#2}{\i};
			\chordPointDSet{b}{#2}{\j};
			\ifthenelse{\i=\j}{
				\coordinate (c) at ($(0,0)!1.25!(a)$);
				\pgfmathsetmacro\chordRHalf{0.25*\chordR}
				\draw[thick, black, l] (c) circle (\chordRHalf);
			}{
				\pgfmathtruncatemacro{\chordaout}{180+90-360*\i/#2}
				\pgfmathtruncatemacro{\chordain}{180+90-360*\j/#2}
				\draw[ultra thick, black, l] (a) to[out=\chordaout,in=\chordain] (b);
			}
		}
		\foreach \i/\j in {#3} {
			\chordPointDSet{a}{#2}{\i};
			\chordPointDSet{b}{#2}{\j};
			\draw[black,c=4pt]
				(a) node[fill] {}
				(b) node[fill] {};
		}
	\end{tikzpicture}%
}
\DeclareRobustCommand{\chordp}[3]{%
    \tikzset{external/export next=false}%
	\begin{tikzpicture}[scale=#1,tcenter]
		\draw[thick, black, l, dotted] \chordCircle;
		\foreach \i/\j in {#3} {
			\chordPointDSet{a}{#2}{\i};
			\chordPointDSet{b}{#2}{\j};
			\ifthenelse{\i=\j}{
				\coordinate (c) at ($(0,0)!1.25!(a)$);
				\pgfmathsetmacro\chordRHalf{0.25*\chordR}
				\draw[thick, black, l] (c) circle (\chordRHalf);
			}{
				\pgfmathtruncatemacro{\chordaout}{180+90-360*\i/#2}
				\pgfmathtruncatemacro{\chordain}{180+90-360*\j/#2}
				\draw[ultra thick, black, l] (a) to[out=\chordaout,in=\chordain] (b);
			}
		}
		\foreach \i/\j in {#3} {
			\chordPointDSet{a}{#2}{\i};
			\chordPointDSet{b}{#2}{\j};
			\draw[black,c=4pt]
				(a) node[fill] {}
				(b) node[fill] {};

			\pgfmathsetmacro\afrom{90-360*\i/#2-270/#2}
			\pgfmathsetmacro\ato{90-360*\i/#2+270/#2}
			\draw[black,thick,l]
				($\chordR * cos(\afrom) *(1,0) + \chordR * sin(\afrom) *(0,1)$) arc (\afrom:\ato:\chordR);

			\pgfmathsetmacro\bfrom{90-360*\j/#2-270/#2}
			\pgfmathsetmacro\bto{90-360*\j/#2+270/#2}
			\draw[black,thick,l]
				($\chordR * cos(\bfrom) *(1,0) + \chordR * sin(\bfrom) *(0,1)$) arc (\bfrom:\bto:\chordR);
		}
	\end{tikzpicture}%
}
\newcommand{\tchord}[2]{\chord{0.4}{#1}{#2}}
\newcommand{\tchordp}[2]{\chordp{0.4}{#1}{#2}}
\def\graphPointSet#1#2#3{\coordinate (#1) at ($#2 *(1,0) + #3 *(0,1)$)}
\DeclareRobustCommand{\graph}[3]{%
    \tikzset{external/export next=false}%
	\begin{tikzpicture}[scale=#1,tcenter]
		\path (0,0) (0,1);
		\foreach \i/\j/\l in {#2} {
			\graphPointSet{\i}{\j}{\l};
			\draw[black,c=4pt]
				(\i) node[fill] {};
		}
		\foreach \i/\j in {#3} {
			\draw[ultra thick, black, l]
				(\i) -- (\j);
		}
	\end{tikzpicture}%
}
\newcommand{\tgraph}[2]{\graph{0.4}{#1}{#2}}

\newcommand{\tbridge}{%
	\kern-0.2ex%
    \tikzset{external/export next=false}%
	\begin{tikzpicture}[scale=0.4,tcenter,use Hobby shortcut]
		\path (0,-1) (0,1);
		\draw[black,very thick,l]
			(-0.3,0.2) .. (0,0.1) .. (0.3,0.2)
			(-0.3,-0.2) .. (0,-0.1) .. (0.3,-0.2);
	\end{tikzpicture}%
	\kern-0.2ex%
}

\begin{document}
\title{The sl(2)-weight system at c~=~3/8 for graphs}
\author[1]{Daniil Fomichev}
\affiliation[1]{
organization = {Saint-Petersburg State University},
postcode = {199034},
addressline ={ 7/9 Universitetskaya Emb.},
city = { Saint-Petersburg},
country ={ Russia}
}
\ead{fomichev.d.s@yandex.ru}
\author[2]{Maksim Karev}
\affiliation[2]{
organization = {Guangdong Technion-Israel Institute of Technology},
postcode = {515603},
addressline = {241 Daxue Lu},
city = {Shantou city, Guangdong province},
country = {P.R. China}
}
\ead{maksim.karev@gtiit.edu.cn}

\date{\today}

\begin{abstract}
We construct a 4-invariant that extends the specialization of the $\mathfrak{sl}(2)$-weight system at $c = \frac 38$ and satisfies a simple deletion-contraction relation.
\end{abstract}


\maketitle

\tableofcontents

\section{Introduction}
\subsection{Motivation and history of the problem}

Weight systems are functions defined on chord diagrams that naturally emerge in the study of finite type invariants of knots (see~\cite{vassiliev1990cohomology} and~\cite{kontsevich1993vassiliev}). 

Metrized Lie (super)algebras serve as one of the primary sources for weight systems (see, for example, ~\cite{kontsevich1993vassiliev, BN95, CDBook}. According to the theorem by S. Chmutov and S. Lando (\cite{ChL07}), the values of weight systems associated with Lie algebra $\mathfrak{sl}(2)$ and Lie superalgebra $\mathfrak{gl}(1|1)$ are determined by the intersection graphs of chord diagrams. The graph-theoretic counterparts of weight systems  are called 4-invariants  (see~\cite{landohopf}), and every 4-invariant 
gives rise
 %corresponds 
to a weight system 
by means of
 %through 
the intersection graph map. The map from the space of 4-invariants to the set of weight systems is 
not surjective. 
 %non-surjective. 
For instance, it is not even clear if a weight system whose values only depend on the intersection graph of the chord diagram
 %with values depending on intersection graphs only, 
can always be extended to a 4-invariant. For the case of $\mathfrak{gl}(1|1)$, an explicit 4-invariant that extends
 %extending 
it is known. A long-standing conjecture of S. Lando asserts:
\begin{hypothesis}
There exists a unique 4-invariant that extends %the value of 
 the $\mathfrak{sl}(2)$-weight system.
\end{hypothesis}

The current state of knowledge %research 
 on the existence of a 4-invariant extending the $\mathfrak{sl}(2)$-weight system can be found in~\cite{KL23} and references therein. We mention the following results: there exists a unique 4-invariant that  extends the $\mathfrak{sl}(2)$-weight system   %and is defined non-trivially on 
 to the space of all graphs with up to eight vertices~\cite{K21}. Some coefficients of the $\mathfrak{sl}(2)$-weight system are related to the numerical invariants of the corresponding intersection graphs~\cite{CDBook,KLMR14,BNV15}. There are also studies discussing the values of possible extensions on particular families of graphs~\cite{Z20,Z22,Za23,ZZ24}. 

While we do not resolve the Lando conjecture, we present an explicit construction of a 4-invariant extending the %a 
 specialization of the $\mathfrak{sl}(2)$-weight system at $c = \frac 38$. %Unfortunately, 
  The question of uniqueness of such an extension is still outside our reach.

The paper is structured as follows: Section 2 provides essential information about the $\mathfrak{sl}(2)$-weight system and 4-invariants. We outline the defining relations for the $\mathfrak{sl}(2)$-weight system and introduce a certain set of relations on
 functions defined over graphs with a view to construct the desired %facilitate an      
 extension. 
 Section 3 details the definition of function $\phi$, the deletion-contraction relation it satisfies, and proofs establishing that this function is a 4-invariant extending the $\mathfrak{sl}(2)$-weight system at $c = \frac 38$. Section 4 contains two more formulae for function $\phi$ and discusses an upper bound on its value.

The key results of this paper include theorem~\ref{thm:4inv} confirming that $\phi$ is a 4-invariant and theorem~\ref{thm: main} establishing that $\phi$ extends the specialization of the $\mathfrak{sl}(2)$-weight system.

After the first version of this text was completed, P.~Zakorko brought to our attention the relation of the constructed invariant $\phi$ to the standard two dimensional representation of $\sl(2)$. We discuss this relation further in remark~\ref{remark:sl2-2d}, which concludes Section 3.

\subsection{Graph drawing conventions}\label{section:graph-convention}

In this paper, we only consider simple graphs. The set of isomorphism classes of simple graphs on $n$ vertices is denoted $\mathbf{G}_n$, and $\mathbf{G} = \bigcup_{n\ge 0} \mathbf{G}_n$.

We follow a particular convention for %visually 
 depicting graphs. Some of the vertices of a graph will be shown individually, others may be omitted, while yet others will be grouped into (not necessarily disjoint) sets that will be shown as circular or elliptic nodes. We will typically denote these sets of vertices by lowercase Latin letters that will be shown inside the nodes.
 If a vertex $v$ is shown to be connected with an edge to a node labeled with $a$, the graph has edges connecting $v$ to all the vertices within the subset $a$ of the vertices of the graph. If a vertex $v$ is shown in the figure, all the edges emanating from $v$ are also shown, either explicitly or as edges connecting $v$ to the nodes.
 %  Not individually shown, certain vertices may be grouped into sets typically represented by lowercase Latin letters. If a vertex $v$ is visibly connected to a node labeled with a letter $a$, the graph has edges connecting vertex $v$ to all vertices within subset $a$ of vertices of the graph.  All the edges that connect an explicitly shown vertex $v$ to other vertices of the graph are either shown explicitly or made part of a group of edges shown by the connection of $v$ to a node. 
 It is important to note that subsets of vertices adjacent to explicitly shown vertices are not required to have empty intersections in this context, %, but 
although sometimes we will explicitly require certain relations between the sets of vertices. When presenting equations involving multiple graphs, it is understood that all the omitted %non-depicted   
vertices and edges, 
as well as the grouping of not explicitly shown vertices  into sets, are consistent across the graphs in question. 

For example, if we only care about the fragment in a neighborhood of the white vertex of the graph depicted below, provided with the following decomposition of the vertices into individual sets, we abbreviate it as follows.
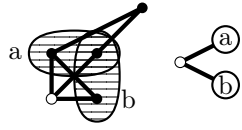
\begin{center}
	\begin{tikzpicture}[scale=0.6,tcenter,use Hobby shortcut]
		\draw[ultra thick, black, l]
			(0,0) -- (1,0)
			(0,0) -- (0,1)
			(0,0) -- (1,1)
			(1,0) -- (0,1)
			(1,1) -- (2,2)
			(0,1) -- (2,2);
		\draw[black,c=4pt]
			(0,0) node[fill=white,draw=black] {}
			(1,0) node[fill] {}
			(0,1) node[fill] {}
			(1,1) node[fill] {}
			(2,2) node[fill] {};
		\draw[thick,l,pattern={Lines[distance=1mm,angle=45,line width=0.1mm]}]
			([closed]-0.5,1) .. (0.5,0.5) .. (1.5,1) .. (0.5,1.5) .. (-0.5,1);
		\draw
			(-0.8,1) node {a};
		\draw[thick,l,pattern={Lines[distance=1mm,angle=-45,line width=0.1mm]}]
			([closed]1,-0.5) .. (0.5,0.5) .. (1,1.5) .. (1.5,0.5) .. (1,-0.5);
		\draw
			(1.7,0) node {b};
	\end{tikzpicture}%
	\quad%
	\begin{tikzpicture}[scale=0.6,tcenter]
		\draw[ultra thick, black, l]
			(1,-0.5) -- (0,0) -- (1,0.5);
		\draw[c=4pt]
			(0,0) node[fill=white,draw=black] {};
		\draw[circle]
			(1,0.5) node[nb] {a};
		\draw[circle]
			(1,-0.5) node[nb] {b};
	\end{tikzpicture}
\end{center}

%On the last picture the colors do not carry any meaningful information and are added for  better comprehension only.

The edges of a
 %the 
 graph are shown as solid. Following~\cite{K24}, we adopt a specific notation for the  values of functions defined on $\mathbf{G}$. Given a function on $\mathbf{G}$, we can extend it by linearity to the vector space $\C \mathbf{G}$, denoting the resulting function %,
   by the same symbol by abuse of notation. We also allow some of the edges of the graph to be dashed, meaning by it the following formal linear combination:  

\begin{align*}
    \begin{tikzpicture}[scale=0.4, tcenter]
        \draw[ultra thick, red pseudo, l] (0, 0) -- (0, 1);
        \draw[black, c=4pt]
            (0, 0) node[fill] {}
            (0, 1) node[fill] {};
    \end{tikzpicture}
    = \begin{tikzpicture}[scale=0.4, tcenter]
        \draw[ultra thick, black, l] (0, 0) -- (0, 1);
        \draw[black, c=4pt]
            (0, 0) node[fill] {}
            (0, 1) node[fill] {};
    \end{tikzpicture}
    - \begin{tikzpicture}[scale=0.4, tcenter]
        \draw[black, c=4pt]
            (0, 0) node[fill] {}
            (0, 1) node[fill] {};
    \end{tikzpicture}.
\end{align*}

When more than one edge is dashed, it denotes the alternating sum over all possible ways of resolving the dashed edges, whether as standard edges or missing edges respectively. A dashed edge connecting an explicitly shown vertex to a node means that all the edges connecting the vertex to the elements of the set %, 
represented by the node %,
are dashed.

\subsection{Acknowledgements}
The authors thank M.E.~Kazarian and J.~Mostovoy for the valuable comments on the first version of this text, E.~Krasilnikov for verification of correctness of the construction of the invariant for all graphs up to 8 vertices, and P. Zakorko for interesting discussions.
We also thank anonymous referees for their careful reading of the manuscript and their insightful comments and suggestions, which have greatly improved the quality of the work.

\section{\texorpdfstring{$\sl(2)$}{sl(2)}-weight system}

This section reviews the necessary definitions and discusses, in brief, the known results about the $\sl(2)$-weight system.

\subsection{Weight systems and 4-invariants}

%The domain of $\sl(2)$-weight system is the set of chord diagrams defined below

\begin{defn}
    For $n\in \N$, a \emph{chord diagram} $D$ of \emph{order} $n$ is a set of pairwise distinct $2n$ points positioned on an oriented circle along with a specific complete pairing of these points. The diagram is considered up to an orientation-preserving diffeomorphism of the circle. We denote the set of all chord diagrams of order $n$ by $\mathbf{A}_n$, and $\mathbf{A} = \bigcup_{n \ge 0} \mathbf{A}_n$.
\end{defn}

In visual representations of chord diagrams, paired points are conventionally connected by a chord, hence the object's name. In this paper, it is assumed that all circles depicted in the figures are oriented counterclockwise.

The V.A. Vassiliev's theory of finite-type knot invariants~\cite{V90, BN95, CDBook} motivates the study of the space of functions defined on the set of chord diagrams that satisfy specific relations known as the 4T-relations.

\begin{defn}\label{rel:4T}
    Let $M$ be an abelian group. A function $f\colon \mathbf{A} \to M$ is called a \emph{weight system} if the following relations hold:

    \begin{align*}
        f\left(\tchordp{16}{0/8,-4/9}\right)
        -f\left(\tchordp{16}{0/8,-4/7}\right)
        +f\left(\tchordp{16}{0/8,-4/1}\right)
        -f\left(\tchordp{16}{0/8,-4/-1}\right)
        =0.
    \end{align*}
\end{defn}

Henceforth, the following convention is used. The diagrams may include additional chords with endpoints on the dashed arcs. However, only the chords with endpoints on the solid parts of the circles are explicitly depicted. All other chords not visible in the images are consistent across all four chord diagrams.

\begin{remark}
    In the theory of non-framed knot invariants, it is also necessary to impose so-called 1T-relations. When 1T-relations are not assumed, the corresponding functions are referred to as ``framed weight systems.'' In this paper, we do not address 1T-relations and simply refer to ``framed weight systems'' as ``weight systems'' for brevity.%textcolor{red}{brevity.}%conciseness.
\end{remark}

\begin{construction}\label{constr:mult}
    Define the product of two chord diagrams as follows: break the supporting circles of each of the two chord diagrams in a point different from chord endpoints and then glue them together, respecting the orientation.

    \begin{align*}
        \tchord{12}{0/8,-2/6,1/5}
        \tbridge
        \tchord{12}{0/4,2/6}
        =\tchord{10}{1/3,2/4,0/5,-1/-3,-2/-4}
        \sim\tchord{10}{2/4,3/5,0/1,-1/-3,-2/-4}
        =\tchord{12}{0/8,-2/6,1/2}
        \tbridge
        \tchord{12}{0/4,2/6}
    \end{align*}

    Different choices of break points, in general, lead to different resulting diagrams. However, the result of the evaluation of a weight system on such a product does not depend on this choice (e.g. see~\cite{CDBook} for details). For chord diagrams $D_1,D_2 \in \mathbf{A}$ and a weight system $f$, by $f(D_1 D_2)$ we denote the result of the evaluation of $f$ on any diagram glued from $D_1$ and $D_2$.
\end{construction}

As it was mentioned in~\cite{chmutovI, chmutovII, chmutov1994vassiliev}, every chord diagram $D$ gives rise to the corresponding \emph{intersection graph} $\G_D$.

\begin{construction}
    The intersection graph $\G_D$ of a chord diagram $D$ is formed as follows: the set of vertices of $V(\G_D)$ is the set of chords of $D$. Two vertices are connected by an edge if the corresponding chords of the chord diagrams intersect (or, equivalently, the corresponding chords endpoints interlace).
\end{construction}
An example of the intersection graph of a chord diagram is shown below.
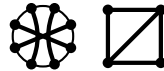
\begin{center}
    \tchord{16}{-1/9,1/7,-3/5,-5/3}%
    \quad \begin{tikzpicture}[scale=0.7, tcenter]
        \draw[ultra thick, black, l]
            (0,0) -- (0,1) -- (1,1) -- (1,0) -- (0,0) -- (1,1);
        \draw[black, c=4pt]
            (0,0) node[fill] {}
            (0,1) node[fill] {}
            (1,0) node[fill] {}
            (1,1) node[fill] {};
    \end{tikzpicture}
\end{center}

\begin{defn}\label{rel:gr4T}
    For an abelian group $M$, a function $g\colon \mathbf{G} \to M$ is called a \emph{4-invariant} if the following  graph 4T-relations hold:
    
    \begin{align*}
        g\left( \begin{tikzpicture}[scale=0.6,tcenter]
            \draw[ultra thick,black,l]
                (0,0) -- (1,0)
                (0,1) -- (1,1);
            \draw[ultra thick,red pseudo,l]
                (0,0) -- (0,1);
            \draw[black,c=4pt]
                (0,0) node[fill] {}
                (0,1) node[fill] {};
            \draw[black,circle]
                (1,1) node[nb] {$a$}
                (1,0) node[nb] {$b$};
        \end{tikzpicture}\right)
        =g\left( \begin{tikzpicture}[scale=0.6,tcenter]
            \draw[ultra thick,black,l]
                (0,0) -- (1.5,0)
                (0,1) -- (1.5,1);
            \draw[ultra thick,red pseudo,l]
                (0,0) -- (0,1);
            \draw[black,c=4pt]
                (0,0) node[fill] {}
                (0,1) node[fill] {};
            \draw[black,ellipse]
                (1.5,1) node[nbe] {$a\triangle b$};
            \draw[black,circle]
                (1.5,0) node[nb] {$b$};
        \end{tikzpicture}\right).
    \end{align*}

    Here we use the graph drawing conventions introduced in subsection~\ref{section:graph-convention}.

	$a\triangle b$ is the symmetric difference of $a$ and $b$.
\end{defn}

\begin{remark}
    Our graph drawing conventions allow us to formulate the graph 4T-relations in a very compact form. We can also state it in a more detailed form, which is closer to the original form  of~\cite{landohopf}:
    \begin{align*}
        g\left( \begin{tikzpicture}[scale=0.6,tcenter]
            \draw[ultra thick,black,l]
                (0,0) -- (1,-0.5)
                (0,1) -- (1,1.5)
                (0,0) -- (1,0.5) -- (0,1);
            \draw[ultra thick,red pseudo,l]
                (0,0) -- (0,1);
            \draw[black,c=4pt]
                (0,0) node[fill] {}
                (0,1) node[fill] {};
            \draw[black,circle]
                (1,1.5) node[nb] {$a$}
                (1,0.5) node[nb] {$c$}
                (1,-0.5) node[nb] {$b$};
        \end{tikzpicture}\right)
        =g\left( \begin{tikzpicture}[scale=0.6,tcenter]
            \draw[ultra thick,black,l]
                (0,0) -- (1,-0.5)
                (0,1) -- (1,1.5)
                (0,0) -- (1,0.5) -- (0,1);
            \draw[ultra thick,red pseudo,l]
                (0,0) -- (0,1);
            \draw[black,c=4pt]
                (0,0) node[fill] {}
                (0,1) node[fill] {};
            \draw[black,circle]
                (1,1.5) node[nb] {$c$}
                (1,0.5) node[nb] {$a$}
                (1,-0.5) node[nb] {$b$};
        \end{tikzpicture}\right).
    \end{align*}
    However, in this case we have to state explicitly that the vertex subsets denoted by $a,b,c$ have pairwise empty intersection.
\end{remark}

Given a 4-invariant $g\colon \mathbf{G} \to M$, we can define a weight system $f$ by the rule: for $D \in \mathbf{A}$, $f(D) = g(\G_D)$. Notice that this construction works due to the fact that application of the intersection graph map to diagrams involved in 4T-relation produces a graph 4T-relation.

This paper deals with weight systems and 4-invariants with values in $\C$ and $\C[t]$ only. For these cases, both the spaces of weight systems and of 4-invariants have a natural bialgebra structure, though we are not going to use it in this paper. For the details see, e.g.,~\cite{CDBook}. We refer to~\cite{KL23} for %as 
 a state-of-the-art review of weight systems and 4-invariants.

\subsection{The Chmutov-Varchenko relations}

According to~\cite{K93, BN95}, every metrized Lie algebra gives rise to a weight system with values in the center of the corresponding universal enveloping algebra. For a particular case of this general construction related to the Lie algebra $\sl(2)$, S. Chmutov and A. Varchenko~\cite{ChV97} have proposed a recursive procedure for computation of the value of the corresponding weight system without any reference to the Lie algebra. The center of the universal enveloping algebra of $\sl(2)$ is isomorphic to the algebra $\C[c]$ of polynomials in one variable, where $c$ is the Casimir element defined by the choice of the metric. Here we choose the metric to be the Killing form of $\sl(2)$\footnote{Any other choice of the metric on $\mathfrak{sl}(2)$ leads to the rescaling of the coefficients of the $\mathfrak{sl}(2)$-weight systems. See~\cite{CDBook} for the details.}. 

\begin{defn}\label{rel:ChV}
	The $\sl(2)$-weight system $w_{\sl(2)}\colon \mathbf{A} \to \C[c]$ is uniquely defined by the following axioms:
    \begin{enumerate}
		\item \emph{Normalization.} The value of $w_{\sl(2)}$ on the chord diagram with one chord equals $c$.
		%\item $w_{\sl(2)}$ is a weight system, i.e. 4T-relations~\ref{rel:4T} hold.
		\item \emph{Multiplicativity.} For two chord diagrams $D_1,D_2 \in \mathbf{A}$, we have $w_{\sl(2)}(D_1 D_2) = w_{\sl(2)}(D_1) w_{\sl(2)}(D_2)$, where $w_{\sl(2)}(D_1 D_2)$ denotes the result of evaluation on any diagram glued from $D_1$ and $D_2$ (see construction~\ref{constr:mult}).
		\item \emph{Leaf deletion.} If chord diagram $D$ has a leaf, i.e. a chord that intersects only one other chord, denoting the result of the deletion of the leaf by $D'$, we have $w_{\sl(2)}(D) = \left(c - \frac{1}{2}\right)w_{\sl(2)}(D')$.
        \item \emph{6T-relations.} The following relations hold:
    
        \begin{align*}
            w_{\sl(2)}\left(\tchordp{22}{0/11,-1/4,12/7}\right)
            - &w_{\sl(2)}\left(\tchordp{22}{0/11,1/4,12/7}\right)
            - w_{\sl(2)}\left(\tchordp{22}{0/11,-1/4,10/7}\right)
            + w_{\sl(2)}\left(\tchordp{22}{0/11,1/4,10/7}\right) = \\
            = &\frac{1}{2} w_{\sl(2)}\left(\tchordp{22}{-1/12,4/7}\right)
            - \frac{1}{2} w_{\sl(2)}\left(\tchordp{22}{-1/7,12/4}\right); \\
            w_{\sl(2)}\left(\tchordp{22}{0/11,-1/7,12/4}\right)
            - &w_{\sl(2)}\left(\tchordp{22}{0/11,1/7,12/4}\right)
            - w_{\sl(2)}\left(\tchordp{22}{0/11,-1/7,10/4}\right)
            + w_{\sl(2)}\left(\tchordp{22}{0/11,1/7,10/4}\right) = \\
            = &\frac{1}{2} w_{\sl(2)}\left(\tchordp{22}{-1/12,4/7}\right)
            - \frac{1}{2} w_{\sl(2)}\left(\tchordp{22}{-1/4,12/7}\right). \\
        \end{align*}
    \end{enumerate}
\end{defn}

The above-listed properties are known as \emph{Chmutov-Varchenko relations} for the $\sl(2)$-weight system. Notably, the 4T-relations are not directly required, but are implied by the defining axioms.

Motivated by the theorem of S. Chmutov and S. Lando~\cite{ChL07}, we introduce the following definition, related to the previous one by translating all identities involving chord diagrams to the intersection graph language.
\begin{defn}[The graph Chmutov-Varchenko  relations]\label{defn:ChV}
    We say that a function $g\colon \mathbf{G} \to \C$ satisfies the graph Chmutov-Varchenko  relations at $p \in \C$ if the following conditions hold:
    \begin{enumerate}
        \item \emph{Normalization.} The value of $g$ on the graph with only one vertex and no edges equals $p$.

        \item \emph{Multiplicativity.} For any two graphs $G_1,G_2 \in \mathbf{G}$, the value of $g$ on their disjoint union is: $g(G_1 \sqcup G_2) = g(G_1) g(G_2)$.
        \item \emph{Leaf deletion.} If $G$ is a graph with a leaf, i.e. a vertex of degree 1, and $G'$ is the result of removing this vertex together with the edge adjacent to it, then $g(G) = \left(p-\frac{1}{2}\right) g(G').$
        \item \emph{Graph 6T-relations.} The following relations hold:
            \begin{align*}
                g\left( \begin{tikzpicture}[scale=0.4,tcenter]
                    \draw[ultra thick, red pseudo, l]
                        (0,1) -- (0,-1);
                    \draw[ultra thick, black, l]
                        (0,1) -- (1.5,1)
                        (0,0) -- (1.5,0)
                        (0,-1) -- (1.5,-1);
                    \draw[black, c=4pt]
                        (0,1) node[fill] {}
                        (0,0) node[fill] {}
                        (0,-1) node[fill] {};
                    \draw[black,circle]
                        (1.5,1) node[nb] {$x$}
                        (1.5,0) node[nb] {$y$}
                        (1.5,-1) node[nb] {$z$};
                \end{tikzpicture}\right)
                &=\frac{1}{2} g\left( \begin{tikzpicture}[scale=0.4,tcenter]
                    \draw[ultra thick, black, l]
                        (0,1) -- (3,1)
                        (3,-1) -- (0,-1);
                    \draw[black, c=4pt]
                        (0,1) node[fill] {}
                        (0,-1) node[fill] {};
                    \draw[black,ellipse]
                        (3,1) node[nb] {$y$}
                        (3,-1) node[nbe=33pt] {$x\triangle y\triangle z$};
                \end{tikzpicture}\right)
                -\frac{1}{2} g\left( \begin{tikzpicture}[scale=0.4,tcenter]
                    \draw[ultra thick, black, l]
                        (0,1) -- (1.5,1)
                        (1.5,-1) -- (0,-1)
                        (0,1) -- (0,-1);
                    \draw[black, c=4pt]
                        (0,1) node[fill] {}
                        (0,-1) node[fill] {};
                    \draw[black,ellipse]
                        (2,1) node[nbe] {$x\triangle y$}
                        (2,-1) node[nbe] {$y\triangle z$};
                \end{tikzpicture}\right); \\
                g\left( \begin{tikzpicture}[scale=0.4,tcenter]
                    \draw[ultra thick, black, l]
                        (0,1) -- (1.5,1)
                        (0,0) -- (1.5,0)
                        (0,-1) -- (1.5,-1)
                        (0,1) edge[out=-120,in=120] (0,-1);
                    \draw[ultra thick, red pseudo, l]
                        (0,1) -- (0,-1);
                    \draw[black, c=4pt]
                        (0,1) node[fill] {}
                        (0,0) node[fill] {}
                        (0,-1) node[fill] {};
                    \draw[black,circle]
                        (1.5,1) node[nb] {$x$}
                        (1.5,0) node[nb] {$y$}
                        (1.5,-1) node[nb] {$z$};
                \end{tikzpicture}\right)
                &=\frac{1}{2} g\left( \begin{tikzpicture}[scale=0.4,tcenter]
                    \draw[ultra thick, black, l]
                        (0,1) -- (3,1)
                        (3,-1) -- (0,-1);
                    \draw[black, c=4pt]
                        (0,1) node[fill] {}
                        (0,-1) node[fill] {};
                    \draw[black,ellipse]
                        (3,1) node[nb] {$y$}
                        (3,-1) node[nbe=33pt] {$x\triangle y\triangle z$};
                \end{tikzpicture}\right)
                -\frac{1}{2} g\left( \begin{tikzpicture}[scale=0.4,tcenter]
                    \draw[ultra thick, black, l]
                        (0,1) -- (1.5,1)
                        (1.5,-1) -- (0,-1);
                    \draw[black, c=4pt]
                        (0,1) node[fill] {}
                        (0,-1) node[fill] {};
                    \draw[black,ellipse]
                        (2,1) node[nbe] {$x\triangle y$}
                        (2,-1) node[nbe] {$y\triangle z$};
                \end{tikzpicture}\right).
            \end{align*}
        %\item $g\colon \mathbf{G} \to \C$ is a 4-invariant, i.e. the graph 4T-relations~\ref{rel:gr4T} hold.
    \end{enumerate}
\end{defn}

\begin{remark}
    A reader may find it not straightforward to relate the 6T-relations from the definition of the $\sl(2)$-weight system to the graph 6T-relations we introduce here. To make this task easier, we give another form of graph 6T-relations assuming that the vertices sets $x,y,z,t,u,v,w$ have pairwise empty intersections:
    \begin{align*}
        g\left( \begin{tikzpicture}[scale=0.4,tcenter]
            \draw[ultra thick, red pseudo, l]
                (0,1) -- (0,-1);
            \draw[ultra thick, black, l]
                (0,0) -- (1,0)
                (0,1) -- (1,0)
                (0,-1) -- (1,0)
                (0,2) -- (0,1) -- (-1,0.5) -- (0,0) -- (-1,-0.5) -- (0,-1) -- (0,-2)
                (-2.5,0) edge[in=150,out=90] (0,1)
                (-2.5,0) edge[in=-150,out=-90] (0,-1)
                (0,0) edge[bend left=60] (2.5,0);
            \draw[black, c=4pt]
                (0,1) node[fill] {}
                (0,0) node[fill] {}
                (0,-1) node[fill] {};
            \draw[black,circle]
                (0,2) node[nb] {$x$}
                (0,-2) node[nb] {$z$}
                (1,0) node[nb] {$t$}
                (2.5,0) node[nb] {$y$}
                (-1,0.5) node[nb] {$v$}
                (-1,-0.5) node[nb] {$w$}
                (-2.5,0) node[nb] {$u$};
        \end{tikzpicture}\right)
        &=\frac{1}{2} g\left( \begin{tikzpicture}[scale=0.4,tcenter]
            \draw[ultra thick, black, l]
                (0,2) -- (0,1) -- (2,0) -- (0,-1) -- (0,-2);
            \draw[black, c=4pt]
                (0,1) node[fill] {}
                (0,-1) node[fill] {};
            \draw[black,circle]
                (-1,0) node[nb] {$u$};
            \draw[black,ellipse]
                (2,0) node[nbe] {$y\kern-0.4ex\cup\kern-0.4ex t$}
                (0,2) node[nbe] {$v\kern-0.4ex\cup\kern-0.4ex w$}
                (0,-2) node[nbe] {$x\kern-0.4ex\cup\kern-0.4ex z$};
        \end{tikzpicture}\right)
        -\frac{1}{2} g\left( \begin{tikzpicture}[scale=0.4,tcenter]
            \draw[ultra thick, black, l]
                (0,2) -- (0,1) -- (2,0) -- (0,-1) -- (0,-2)
                (0,1) -- (0,-1);
            \draw[black, c=4pt]
                (0,1) node[fill] {}
                (0,-1) node[fill] {};
            \draw[black,circle]
                (-1,0) node[nb] {$t$};
            \draw[black,ellipse]
                (0,2) node[nbe] {$v\kern-0.4ex\cup\kern-0.4ex z$}
                (0,-2) node[nbe] {$x\kern-0.4ex\cup\kern-0.4ex w$}
                (2,0) node[nbe] {$u\kern-0.4ex\cup\kern-0.4ex y$};
        \end{tikzpicture}\right); \\
        g\left( \begin{tikzpicture}[scale=0.4,tcenter]
            \draw[ultra thick, red pseudo, l]
                (0,1) -- (0,-1);
            \draw[ultra thick, black, l]
                (0,1) -- (1,0) -- (0,-1)
                (0,0) -- (1,0)
                (0,2) -- (0,1) -- (-1,0.5) -- (0,0) -- (-1,-0.5) -- (0,-1) -- (0,-2)
                (-2.5,0) edge[in=150,out=90] (0,1)
                (-2.5,0) edge[in=-150,out=-90] (0,-1)
                (0,1) edge[out=-60,in=60] (0,-1)
                (0,0) edge[bend left=60] (2.5,0);
            \draw[black, c=4pt]
                (0,1) node[fill] {}
                (0,0) node[fill] {}
                (0,-1) node[fill] {};
            \draw[black,circle]
                (0,2) node[nb] {$x$}
                (2.5,0) node[nb] {$y$}
                (0,-2) node[nb] {$z$}
                (1,0) node[nb] {$t$}
                (-1,0.5) node[nb] {$v$}
                (-1,-0.5) node[nb] {$w$}
                (-2.5,0) node[nb] {$u$};
        \end{tikzpicture}\right)
        &=\frac{1}{2} g\left( \begin{tikzpicture}[scale=0.4,tcenter]
            \draw[ultra thick, black, l]
                (0,2) -- (0,1) -- (2,0) -- (0,-1) -- (0,-2);
            \draw[black, c=4pt]
                (0,1) node[fill] {}
                (0,-1) node[fill] {};
            \draw[black,circle]
                (-1,0) node[nb] {$u$};
            \draw[black,ellipse]
                (2,0) node[nbe] {$y\kern-0.4ex\cup\kern-0.4ex t$}
                (0,2) node[nbe] {$v\kern-0.4ex\cup\kern-0.4ex w$}
                (0,-2) node[nbe] {$x\kern-0.4ex\cup\kern-0.4ex z$};
        \end{tikzpicture}\right)
        -\frac{1}{2} g\left( \begin{tikzpicture}[scale=0.4,tcenter]
            \draw[ultra thick, black, l]
                (0,2) -- (0,1) -- (2,0) -- (0,-1) -- (0,-2);
            \draw[black, c=4pt]
                (0,1) node[fill] {}
                (0,-1) node[fill] {};
            \draw[black,circle]
                (-1,0) node[nb] {$t$};
            \draw[black,ellipse]
                (2,0) node[nbe] {$u\kern-0.4ex\cup\kern-0.4ex y$}
                (0,2) node[nbe] {$v\kern-0.4ex\cup\kern-0.4ex z$}
                (0,-2) node[nbe] {$x\kern-0.4ex\cup\kern-0.4ex w$};
        \end{tikzpicture}\right).
    \end{align*}
    
    Compared to~definition~\ref{rel:ChV}, for the graph Chmutov-Varchenko relations we do not assume $t=\emptyset$ in the first relation and $y=\emptyset$ in the second.
 
    The graph Chmutov-Varchenko relations take their origin from translating the Chmutov-Varchenko relations for chord diagrams to graphs via the intersection graph map and generalize them: every function that satisfies the graph Chmutov-Varchenko relations defines an extension of the specialization of the $\sl(2)$-weight system at the point $c = p$.
\end{remark}

\section{Definition and principal properties of function \texorpdfstring{$\phi$}{ φ }}
This section discusses the main protagonist of this paper. We introduce the function $\phi\colon \mathbf{G} \to \C$ and prove that it satisfies the graph Chmutov-Varchenko relations (see definition~\ref{defn:ChV}) at $\frac{3}{8}$.

\subsection{Definition of function \texorpdfstring{$\phi$}{ φ }}

For graph $G\in \mathbf{G}$, we denote the set of its vertices by $V(G)$ and the set of its edges by $E(G)$. For subset $E' \subset E(G)$, we denote by $G|_{E'}$ the spanning subgraph of $G$ with the set of edges $E'$. By $\chi_3(G)$ we mean the number of proper colorings of the vertices of $G$ in three colors. Recall that a function $V(G)\to \setp{1,\ldots,k}$ is a (proper) coloring of graph $G$ if no two vertices connected by an edge are assigned the same value.

\begin{defn}[Function $\phi$]\label{defn:phi}
    For $G\in \mathbf{G}$, we set
    $$\phi(G) = 2^{-3|V(G)|}\sum_{E' \subset E(G)}(-2)^{|E'|} \chi_3(G|_{E'}).$$
\end{defn}

Here are the values of $\phi$ on some simple graphs:
\begin{itemize}
    \item For the discrete graph $N_n$ on $n$ vertices, we have $\phi(N_n) = \left(\frac{3}{8}\right)^n$.
    \item $\phi\left(\tgraph{1/0/0,2/0/1}{1/2}\right) = -\frac{3}{8^2}$.
    \item $\phi\left(\tgraph{1/0/0,2/1/0.5,3/0/1}{1/2,2/3}\right) = \frac{3}{8^3}$.
    \item $\phi\left(\tgraph{1/0/0,2/1/0.5,3/0/1}{1/2,2/3,3/1}\right) = \frac{15}{8^3}$.
\end{itemize}

Note that if we compute the value of $\phi$ on graph $G$ with the subset $R$ of edges being dashed (following the conventions from subsection~\ref{section:graph-convention}), we obtain:
$$\phi(G) = 2^{-3|V(G)|}\sum_{E':\ R \subset E' \subset E(G)}(-2)^{|E'|} \chi_3(G|_{E'}).$$
Here the notion of proper vertex colorings of the graph does not depend on the type of edges connecting the vertices. The formula reduces to a single term and becomes especially simple when we evaluate the function on the graph with dashed edges only.

The remaining part of this section is dedicated to the proof of the fact that $\phi$ satisfies the graph Chmutov-Varchenko  relations at  $\frac{3}{8}$. Two out of four properties are straightforward.

\begin{prop}\label{prop:mulnorm}
    Function $\phi$ is multiplicative and is equal to $\frac{3}{8}$ on the graph with one vertex and no edges.
\end{prop}
\begin{proof}
    The value of $\phi$ on the graph with one vertex and no edges is trivially computed by definition.
    Multiplicativity follows from the fact that for a disjoint union $G = G_1 \sqcup G_2$, the set of subsets of edges of $G$ is canonically isomorphic to the direct product of the sets of subsets of edges of $G_1$ and $G_2$.
    Furthermore, the number of vertices and the number of edges of a graph are both additive under the disjoint union, while the number of proper 3-colorings is multiplicative. The assertion follows.
\end{proof}

\subsection{The deletion-contraction identity for \texorpdfstring{$\phi$}{ φ }}

Function $\phi$ takes its origin in the number of proper vertex colorings of the graph in three colors. Unsurprisingly, a version of the deletion-contraction relation allowing one to compute its value recursively can be stated.

\begin{thm}\label{thm:del-cont}
    The following deletion-contraction relation holds:

    \begin{align}\label{rel:del-cont}
        \phi\left( \begin{tikzpicture}[scale=0.6,tcenter]
            \draw[ultra thick,black,l]
                (1,1) -- (0,1) -- (0,0) -- (1,0);
            \draw[black,c=4pt]
                (0,0) node[fill] {}
                (0,1) node[fill] {};
            \draw[black,circle]
                (1,1) node[nb] {$u$}
                (1,0) node[nb] {$v$};
        \end{tikzpicture}\right)
        =-\phi\left( \begin{tikzpicture}[scale=0.6,tcenter]
            \draw[ultra thick,black,l]
                (1,1) -- (0,1) (0,0) -- (1,0);
            \draw[black,c=4pt]
                (0,0) node[fill] {}
                (0,1) node[fill] {};
            \draw[black,circle]
                (1,1) node[nb] {$u$}
                (1,0) node[nb] {$v$};
        \end{tikzpicture}\right)
        +\frac{1}{4} \phi\left( \begin{tikzpicture}[scale=0.6,tcenter]
            \draw[ultra thick,black,l]
                (0,0) -- (1.5,0);
            \draw[black,c=4pt]
                (0,0) node[fill] {};
            \draw[black,ellipse]
                (1.5,0) node[nbe] {$u\triangle v$};
        \end{tikzpicture}\right).
    \end{align}
\end{thm}

\begin{proof}
    Recall the definition of $\phi$:
    \begin{align*}
        \phi(G)
        &= 2^{-3|V(G)|} \sum_{E'\subset E(G)} (-2)^{|E'|} \chi_3(G|_{E'})
    \end{align*}
    To specify every term of the sum separately, we can write down the graph under consideration as a linear combination of graphs in which some of the edges are dashed:
    \begin{align*}
        \phi\left( \begin{tikzpicture}[scale=0.6,tcenter]
            \draw[ultra thick,black,l]
                (1,1) -- (0,1) -- (0,0) -- (1,0);
            \draw[black,c=4pt]
                (0,0) node[fill] {}
                (0,1) node[fill] {};
            \draw[black,circle]
                (1,1) node[nb] {$u$}
                (1,0) node[nb] {$v$};
        \end{tikzpicture}\right) = 
        \sum_{\begin{smallmatrix}u' \subset u \\ v' \subset v \end{smallmatrix}}
        \phi\left( \begin{tikzpicture}[scale=0.6,tcenter]
            \draw[ultra thick,red pseudo,l]
                (1,1) -- (0,1) -- (0,0) -- (1,0);
            \draw[black,c=4pt]
                (0,0) node[fill] {}
                (0,1) node[fill] {};
            \draw[black,circle]
                (1,1) node[nb] {$u'$}
                (1,0) node[nb] {$v'$};
        \end{tikzpicture}\right)
        + \sum_{\begin{smallmatrix}u' \subset u \\ v' \subset v \end{smallmatrix}}
        \phi\left( \begin{tikzpicture}[scale=0.6,tcenter]
            \draw[ultra thick,red pseudo,l]
                (1,1) -- (0,1);
            \draw[ultra thick,red pseudo,l]
                (0,0) -- (1,0);
            \draw[black,c=4pt]
                (0,0) node[fill] {}
                (0,1) node[fill] {};
            \draw[black,circle]
                (1,1) node[nb] {$u'$}
                (1,0) node[nb] {$v'$};
        \end{tikzpicture}\right).
    \end{align*}
    We use Tutte's deletion-contraction relation for the chromatic polynomial (see, e.g.,~\cite{dong2005chromatic} and also~\cite{Tutte} for the history of its discovery) applied to the edge between two explicitly shown vertices to compute the contribution from the first term on the right-hand side of the equation. The deletion term and the second term of the right-hand side of the equation produce together
     \begin{align*}
        -\sum_{\begin{smallmatrix}u' \subset u \\ v' \subset v \end{smallmatrix}}
        \phi\left( \begin{tikzpicture}[scale=0.6,tcenter]
            \draw[ultra thick,red pseudo,l]
                (1,1) -- (0,1);
            \draw[ultra thick,red pseudo,l]
                (0,0) -- (1,0);
            \draw[black,c=4pt]
                (0,0) node[fill] {}
                (0,1) node[fill] {};
            \draw[black,circle]
                (1,1) node[nb] {$u'$}
                (1,0) node[nb] {$v'$};
        \end{tikzpicture}\right),
    \end{align*}
    as the deletion term is provided with the factor of $(-2)$ since after the deletion the number of edges decreases by one.
    
    For the contraction term we have the following: consider a vertex $a \in u'\cap v'$. Notice that the results of contraction of the terms corresponding to the pairs of subsets $(u',v')$, $(u'\setminus \{a\}, v'),$ and $(u', v'\setminus \{a\})$ all produce the same graphs, but the number of edges in the graph that corresponds to $(u',v')$ is one greater than for the other two pairs of subsets, so it is counted with an additional factor of $-2$. It means that we can only have a non-vanishing contribution from the terms with $u' \subset u \setminus v$ and $v' \subset v \setminus u$, or, equivalently, only the subsets $w = u'\cup v' \subset u\triangle v$ contribute non-trivially to the contraction term. After the contraction of the edge, the resulting graph has one edge and one vertex less, which leads to the overall factor of $\frac {-2}8 = - \frac 14$. The Theorem follows by summing up all the contributions and rewriting the graphs with dashed edges as linear combinations of graphs with solid edges.

\end{proof}

Although it is not difficult to prove the following fact directly from the definition of $\phi$, as a corollary of deletion-contraction we have:
\begin{prop}\label{prop:ldel}
    The invariant $\phi$ satisfies the leaf deletion property of~\ref{defn:ChV} for $p = \frac{3}{8}$.
\end{prop}
\begin{proof}
    Apply the deletion-contraction relation~(\ref{rel:del-cont}) to the leaf. The deletion term develops an isolated vertex that, by the multiplicativity property, can be replaced by the factor $\frac{3}{8}$. The resulting graph coincides with the contraction term. Thus the total contribution equals $-\frac{3}{8} + \frac{1}{4} = -\frac{1}{8} = \frac{3}{8} - \frac{1}{2}$. The assertion follows.
\end{proof}

\begin{lemma}
    The following identity  holds:
    
    \begin{align}\label{rel:del-cont-var}
        \phi\left( \begin{tikzpicture}[scale=0.6,tcenter]
            \draw[ultra thick,red pseudo,l]
                (0,1) -- (0.5,0.5) -- (0,0);
            \draw[ultra thick,black,l]
                (1,1) -- (0,1) -- (0,0) -- (1,0)
                (0.5,0.5) -- (1.5,0.5);
            \draw[black,c=4pt]
                (0,0) node[fill] {}
                (0,1) node[fill] {}
                (0.5,0.5) node[fill] {};
            \draw[black,circle]
                (1,1) node[nb] {$u$}
                (1,0) node[nb] {$v$}
                (1.5,0.5) node[nb] {$w$};
        \end{tikzpicture}\right)
        =-\phi\left( \begin{tikzpicture}[scale=0.6,tcenter]
            \draw[ultra thick,red pseudo,l]
                (0,1) -- (0.5,0.5) -- (0,0);
            \draw[ultra thick,black,l]
                (1,1) -- (0,1) (0,0) -- (1,0)
                (0.5,0.5) -- (1.5,0.5);
            \draw[black,c=4pt]
                (0,0) node[fill] {}
                (0,1) node[fill] {}
                (0.5,0.5) node[fill] {};
            \draw[black,circle]
                (1,1) node[nb] {$u$}
                (1,0) node[nb] {$v$}
                (1.5,0.5) node[nb] {$w$};
        \end{tikzpicture}\right)
        -\frac{1}{2} \phi\left( \begin{tikzpicture}[scale=0.6,tcenter]
            \draw[ultra thick,red pseudo,l]
                (0,0) -- (0,1);
            \draw[ultra thick,black,l]
                (0,0) -- (1.5,0)
                (0,1) -- (1.5,1);
            \draw[black,c=4pt]
                (0,0) node[fill] {}
                (0,1) node[fill] {};
            \draw[black,circle]
                (1.5,1) node[nb] {$w$};
            \draw[black,ellipse]
                (1.5,0) node[nbe] {$u\triangle v$};
        \end{tikzpicture}\right).
    \end{align}
\end{lemma}
\begin{proof}
    Expressing the graphs with dashed edges in terms of graphs with solid edges by the convention of subsection~\ref{section:graph-convention} and applying the deletion-contraction relation~(\ref{rel:del-cont}) to the resulting graphs, we derive the equality.
\end{proof}

\subsection{\texorpdfstring{$\phi$}{φ} as a 4-invariant}

This subsection discusses the proof of the 4-invariance of $\phi$. In fact,~\cite{K24} contains a proof of the 4-invariance of a function that differs from $\phi$ by a controllable factor, but we find it more instructive to provide an alternative proof here.  We  start with the following lemma:

\begin{lemma}[Triangle identity]
    The following identity holds:
    \begin{align}\label{rel:triangle}
        \phi\left( \begin{tikzpicture}[scale=0.6,tcenter]
            \draw[ultra thick,black,l,red pseudo]
                (-0.5,0) -- (0.5,0) -- (0,0.8) -- (-0.5,0);
            \draw[ultra thick,black,l]
                (-0.5,0) -- (-1,-0.5)
                (0.5,0) -- (1,-0.5)
                (0,0.8) -- (0,1.6);
            \draw[black,c=4pt]
                (-0.5,0) node[fill] {}
                (0,0.8) node[fill] {}
                (0.5,0) node[fill] {};
            \draw[black,circle]
                (-1,-0.5) node[nb] {$x$}
                (0,1.6) node[nb] {$y$}
                (1,-0.5) node[nb] {$z$};
        \end{tikzpicture}\right)
        =\frac{1}{2}\phi\left( \begin{tikzpicture}[scale=0.6,tcenter]
            \draw[ultra thick,black,l,red pseudo]
                (0,0) -- (0,1);
            \draw[ultra thick,black,l]
                (0,0) -- (1,0)
                (0,1) -- (1,1);
            \draw[black,c=4pt]
                (0,0) node[fill] {}
                (0,1) node[fill] {};
            \draw[black,ellipse]
                (1.5,1) node[nbe] {$x\triangle y$}
				(1.5,0) node[nbe] {$y\triangle z$};
		\end{tikzpicture}\right).
    \end{align}
\end{lemma}

\begin{proof}

We prove this identity by induction on the size of $y$. It is enough to consider the case when vertices from $y$ are connected to the explicitly shown vertices of the triangle on the left-hand side of the equation by dashed edges. The induction base $|y| = 0$:
    \begin{align*}
        \phi\left( \begin{tikzpicture}[scale=0.6,tcenter]
            \draw[ultra thick,black,l,red pseudo]
                (-0.5,0) -- (0.5,0) -- (0,0.8) -- (-0.5,0);
            \draw[ultra thick,black,l]
                (-0.5,0) -- (-1,-0.5)
                (0.5,0) -- (1,-0.5);
            \draw[black,c=4pt]
                (-0.5,0) node[fill] {}
                (0,0.8) node[fill] {}
                (0.5,0) node[fill] {};
            \draw[black,circle]
                (-1,-0.5) node[nb] {$x$}
                (1,-0.5) node[nb] {$z$};
        \end{tikzpicture}\right)
        =\frac{1}{2}\phi\left( \begin{tikzpicture}[scale=0.6,tcenter]
            \draw[ultra thick,black,l,red pseudo]
                (0,0) -- (0,1);
            \draw[ultra thick,black,l]
                (0,0) -- (1,0)
                (0,1) -- (1,1);
            \draw[black,c=4pt]
                (0,0) node[fill] {}
                (0,1) node[fill] {};
            \draw[black,circle]
                (1.0,1) node[nb] {$x$}
				(1.0,0) node[nb] {$z$};
		\end{tikzpicture}\right)
    \end{align*}
follows from the definition of $\phi$: in any proper vertex 3-coloring, the color of the topmost vertex of the graph in the left-hand part of the equality is uniquely determined by the colors of the adjacent vertices, so there is a bijection between the sets of proper 3-colorings of the graphs on the left-hand side and the right-hand side. The factor $\frac 12 = \frac {(-2)^2}8$ arises because the graph on the left has two edges and one vertex more than the one on the right.

By induction hypothesis, we assume that for all graphs with  $|y| < n$ for $n\in \mathbb N$ the weighted sum of proper 3-colorings coincide on the left and the right-hand side of the equation. Take a graph with $|y| = n - 1$ and vertex $s \not \in y$ connected to the topmost vertex of the triangle. On all pictures below vertex $s$ is shown explicitly in the middle of the triangle. It will be convenient for us to treat the case in which all edges that connect $a$ to the vertices of the triangle are dashed. Consider the following cases:
\begin{itemize}
\item If $s$ is connected to the topmost vertex of the triangle only, we have to check that 
\begin{align*}
\phi\left( \begin{tikzpicture}[scale=0.6,tcenter]
            \draw[ultra thick,black,l]
                (0,0.3) -- (1,0.8);
            \draw[ultra thick,black,l,red pseudo]
                (-0.5,0) -- (0.5,0) -- (0,0.8) -- (-0.5,0)
                (0,0.8) -- (0,0.3);
            \draw[ultra thick,black,l]
                (-0.5,0) -- (-1,-0.5)
                (0.5,0) -- (1,-0.5)
                (0,0.8) -- (0,1.6);
            \draw[black,c=4pt]
                (-0.5,0) node[fill] {}
                (0,0.8) node[fill] {}
                (0.5,0) node[fill] {}
                (0,0.3) node[fill] {};
            \draw[black,circle]
                (-1,-0.5) node[nb] {$x$}
                (0,1.6) node[nb] {$y$}
                (1,-0.5) node[nb] {$z$}
                (1,0.8) node[nb] {$w$};
        \end{tikzpicture}\right)
= \phi\left( \begin{tikzpicture}[scale=0.6,tcenter]
            \draw[ultra thick,black,l]
                (0,0.8) -- (0,0.3) -- (1,0.8);
            \draw[ultra thick,black,l,red pseudo]
                (-0.5,0) -- (0.5,0) -- (0,0.8) -- (-0.5,0);
            \draw[ultra thick,black,l]
                (-0.5,0) -- (-1,-0.5)
                (0.5,0) -- (1,-0.5)
                (0,0.8) -- (0,1.6);
            \draw[black,c=4pt]
                (-0.5,0) node[fill] {}
                (0,0.8) node[fill] {}
                (0.5,0) node[fill] {}
                (0,0.3) node[fill] {};
            \draw[black,circle]
                (-1,-0.5) node[nb] {$x$}
                (0,1.6) node[nb] {$y$}
                (1,-0.5) node[nb] {$z$}
                (1,0.8) node[nb] {$w$};
        \end{tikzpicture}\right)
-\phi\left( \begin{tikzpicture}[scale=0.6,tcenter]
            \draw[ultra thick,black,l]
                (0,0.3) -- (1,0.8);
            \draw[ultra thick,black,l,red pseudo]
                (-0.5,0) -- (0.5,0) -- (0,0.8) -- (-0.5,0);
            \draw[ultra thick,black,l]
                (-0.5,0) -- (-1,-0.5)
                (0.5,0) -- (1,-0.5)
                (0,0.8) -- (0,1.6);
            \draw[black,c=4pt]
                (-0.5,0) node[fill] {}
                (0,0.8) node[fill] {}
                (0.5,0) node[fill] {}
                (0,0.3) node[fill] {};
            \draw[black,circle]
                (-1,-0.5) node[nb] {$x$}
                (0,1.6) node[nb] {$y$}
                (1,-0.5) node[nb] {$z$}
                (1,0.8) node[nb] {$w$};
        \end{tikzpicture}\right)
\end{align*}
equals
\begin{align*}
\frac 12 \phi\left( \begin{tikzpicture}[scale=0.6,tcenter]
            \draw[ultra thick,black,l]
                (3,0.5) -- (0.5,0.5) -- (0,1)
                (0.5,0.5) -- (0,0);
            \draw[ultra thick,black,l,red pseudo]
                (0,0) -- (0,1);
            \draw[ultra thick,black,l]
                (0,0) -- (1,0)
                (0,1) -- (1,1);
            \draw[black,c=4pt]
                (0,0) node[fill] {}
                (0,1) node[fill] {}
                (0.5,0.5) node[fill] {};
             \draw[black,circle]
                (3,0.5) node[nb] {$w$};
             \draw[black,ellipse]
                (1.5,1) node[nbe] {$x\triangle y$}
				(1.5,0) node[nbe] {$y\triangle z$};
		\end{tikzpicture}\right)
-\frac 12	\phi\left( \begin{tikzpicture}[scale=0.6,tcenter]
            \draw[ultra thick,black,l]
                (3,0.5) -- (0.5,0.5);
            \draw[ultra thick,black,l,red pseudo]
                (0,0) -- (0,1);
            \draw[ultra thick,black,l]
                (0,0) -- (1,0)
                (0,1) -- (1,1);
            \draw[black,c=4pt]
                (0,0) node[fill] {}
                (0,1) node[fill] {}
                (0.5,0.5) node[fill] {};
            \draw[black,circle]
                (3,0.5) node[nb] {$w$};
             \draw[black,ellipse]
                (1.5,1) node[nbe] {$x\triangle y$}
				(1.5,0) node[nbe] {$y\triangle z$};
		\end{tikzpicture}\right),
\end{align*}
so one of the ways to prove this equality is to show that 
\begin{align*}
\phi\left( \begin{tikzpicture}[scale=0.6,tcenter]
            \draw[ultra thick,black,l]
                (0,0.8) -- (0,0.3) -- (1,0.8);
            \draw[ultra thick,black,l,red pseudo]
                (-0.5,0) -- (0.5,0) -- (0,0.8) -- (-0.5,0);
            \draw[ultra thick,black,l]
                (-0.5,0) -- (-1,-0.5)
                (0.5,0) -- (1,-0.5)
                (0,0.8) -- (0,1.6);
            \draw[black,c=4pt]
                (-0.5,0) node[fill] {}
                (0,0.8) node[fill] {}
                (0.5,0) node[fill] {}
                (0,0.3) node[fill] {};
            \draw[black,circle]
                (-1,-0.5) node[nb] {$x$}
                (0,1.6) node[nb] {$y$}
                (1,-0.5) node[nb] {$z$}
                (1,0.8) node[nb] {$w$};
        \end{tikzpicture}\right)
=\phi\left( \begin{tikzpicture}[scale=0.6,tcenter]
            \draw[ultra thick,black,l]
                (0,0.3) -- (1,0.8)
                (-0.5,0) -- (0,0.3) -- (0.5,0);
            \draw[ultra thick,black,l,red pseudo]
                (-0.5,0) -- (0.5,0) -- (0,0.8) -- (-0.5,0);
            \draw[ultra thick,black,l]
                (-0.5,0) -- (-1,-0.5)
                (0.5,0) -- (1,-0.5)
                (0,0.8) -- (0,1.6);
            \draw[black,c=4pt]
                (-0.5,0) node[fill] {}
                (0,0.8) node[fill] {}
                (0.5,0) node[fill] {}
                (0,0.3) node[fill] {};
            \draw[black,circle]
                (-1,-0.5) node[nb] {$x$}
                (0,1.6) node[nb] {$y$}
                (1,-0.5) node[nb] {$z$}
                (1,0.8) node[nb] {$w$};
        \end{tikzpicture}\right),
\end{align*}
or
\begin{align*}
\phi\left( \begin{tikzpicture}[scale=0.6,tcenter]
            \draw[ultra thick,black,l]
                 (0,0.3) -- (1,0.8);
            \draw[ultra thick,black,l,red pseudo]
                (0,0.3) --  (0,0.8) -- (-0.5,0) -- (0.5,0) -- (0,0.8);
            \draw[ultra thick,black,l]
                (-0.5,0) -- (-1,-0.5)
                (0.5,0) -- (1,-0.5)
                (0,0.8) -- (0,1.6);
            \draw[black,c=4pt]
                (-0.5,0) node[fill] {}
                (0,0.8) node[fill] {}
                (0.5,0) node[fill] {}
                (0,0.3) node[fill] {};
            \draw[black,circle]
                (-1,-0.5) node[nb] {$x$}
                (0,1.6) node[nb] {$y$}
                (1,-0.5) node[nb] {$z$}
                (1,0.8) node[nb] {$w$};
        \end{tikzpicture}\right) = 
\phi\left( \begin{tikzpicture}[scale=0.6,tcenter]
            \draw[ultra thick,black,l]
                (0,0.3) -- (1,0.8);
            \draw[ultra thick,black,l,red pseudo]
                (-0.5,0) -- (0.5,0) -- (0,0.8) -- (-0.5,0) -- (0,0.3) -- (0.5,0);
            \draw[ultra thick,black,l]
                (-0.5,0) -- (-1,-0.5)
                (0.5,0) -- (1,-0.5)
                (0,0.8) -- (0,1.6);
            \draw[black,c=4pt]
                (-0.5,0) node[fill] {}
                (0,0.8) node[fill] {}
                (0.5,0) node[fill] {}
                (0,0.3) node[fill] {};
            \draw[black,circle]
                (-1,-0.5) node[nb] {$x$}
                (0,1.6) node[nb] {$y$}
                (1,-0.5) node[nb] {$z$}
                (1,0.8) node[nb] {$w$};
        \end{tikzpicture}\right)
+ \phi\left( \begin{tikzpicture}[scale=0.6,tcenter]
            \draw[ultra thick,black,l]
                (0,0.3) -- (1,0.8);
            \draw[ultra thick,black,l,red pseudo]
                (-0.5,0) -- (0.5,0) -- (0,0.8) -- (-0.5,0) -- (0,0.3);
            \draw[ultra thick,black,l]
                (-0.5,0) -- (-1,-0.5)
                (0.5,0) -- (1,-0.5)
                (0,0.8) -- (0,1.6);
            \draw[black,c=4pt]
                (-0.5,0) node[fill] {}
                (0,0.8) node[fill] {}
                (0.5,0) node[fill] {}
                (0,0.3) node[fill] {};
            \draw[black,circle]
                (-1,-0.5) node[nb] {$x$}
                (0,1.6) node[nb] {$y$}
                (1,-0.5) node[nb] {$z$}
                (1,0.8) node[nb] {$w$};
        \end{tikzpicture}\right)
+ \phi\left( \begin{tikzpicture}[scale=0.6,tcenter]
            \draw[ultra thick,black,l]
                (0,0.3) -- (1,0.8);
            \draw[ultra thick,black,l,red pseudo]
                (0.5,0) -- (-0.5,0) -- (0,0.8) -- (0.5,0) -- (0,0.3);
            \draw[ultra thick,black,l]
                (-0.5,0) -- (-1,-0.5)
                (0.5,0) -- (1,-0.5)
                (0,0.8) -- (0,1.6);
            \draw[black,c=4pt]
                (-0.5,0) node[fill] {}
                (0,0.8) node[fill] {}
                (0.5,0) node[fill] {}
                (0,0.3) node[fill] {};
            \draw[black,circle]
                (-1,-0.5) node[nb] {$x$}
                (0,1.6) node[nb] {$y$}
                (1,-0.5) node[nb] {$z$}
                (1,0.8) node[nb] {$w$};
        \end{tikzpicture}\right).
\end{align*}
Consider the contribution of possible colorings of vertex $s$ to the value of $\phi$ on both sides of the equation. On the left-hand side, only the coloring in which the color of $s$ differs from the color of the topmost vertex can contribute non-trivially. On the right-hand side, due to the factor of $(-2)^{|E'|}$ in the definition of $\phi$, the colorings in which $s$ is colored to the color of the topmost vertex cancel out. The remaining possible colorings of the right-hand side can be put to a bijection to the possible colorings of the left-hand side, so the equality holds.
\item If $s$ is connected to the topmost and only one of the bottom vertices of the triangle, we may use identity~(\ref{rel:del-cont-var}) to reduce the considerations to the previous case:
\begin{align*}
\phi\left( \begin{tikzpicture}[scale=0.6,tcenter]
            \draw[ultra thick,black,l]
                (0,0.3) -- (1,0.8);
            \draw[ultra thick,black,l,red pseudo]
                (-0.5,0) -- (0.5,0) -- (0,0.8) -- (-0.5,0) -- (0,0.3) -- (0,0.8);
            \draw[ultra thick,black,l]
                (-0.5,0) -- (-1,-0.5)
                (0.5,0) -- (1,-0.5)
                (0,0.8) -- (0,1.6);
            \draw[black,c=4pt]
                (-0.5,0) node[fill] {}
                (0,0.8) node[fill] {}
                (0.5,0) node[fill] {}
                (0,0.3) node[fill] {};
            \draw[black,circle]
                (-1,-0.5) node[nb] {$x$}
                (0,1.6) node[nb] {$y$}
                (1,-0.5) node[nb] {$z$}
                (1,0.8) node[nb] {$w$};
        \end{tikzpicture}\right) = 
        -2\phi\left( \begin{tikzpicture}[scale=0.6,tcenter]
            \draw[ultra thick,black,l]
                (0,0.3) -- (1,0.8);
            \draw[ultra thick,black,l,red pseudo]
                (-0.5,0) -- (0.5,0) -- (0,0.8) -- (-0.5,0)
                (0,0.8) -- (0,0.3);
            \draw[ultra thick,black,l]
                (-0.5,0) -- (-1,-0.5)
                (0.5,0) -- (1,-0.5)
                (0,0.8) -- (0,1.6);
            \draw[black,c=4pt]
                (-0.5,0) node[fill] {}
                (0,0.8) node[fill] {}
                (0.5,0) node[fill] {}
                (0,0.3) node[fill] {};
            \draw[black,circle]
                (-1,-0.5) node[nb] {$x$}
                (0,1.6) node[nb] {$y$}
                (1,-0.5) node[nb] {$z$}
                (1,0.8) node[nb] {$w$};
        \end{tikzpicture}\right)
        -\frac 12 \phi\left( \begin{tikzpicture}[scale=0.6,tcenter]
            \draw[ultra thick,black,l,red pseudo]
                (-0.5,0) -- (0.5,0) -- (0,0.8) -- (-0.5,0);
            \draw[ultra thick,black,l]
                (-0.5,0) -- (-1.5,-0.5)
                (0.5,0) -- (1,-0.5)
                (0,0.8) -- (0,1.6);
            \draw[black,c=4pt]
                (-0.5,0) node[fill] {}
                (0,0.8) node[fill] {}
                (0.5,0) node[fill] {};
            \draw[black,circle]
                (0,1.6) node[nb] {$y$}
                (1,-0.5) node[nb] {$z$};
			\draw[black,ellipse]            
            	 (-1.5,-0.5) node[nbe] {$x\triangle w$};    
        \end{tikzpicture}\right).
\end{align*}

The left-hand part of the indicated equality should be equal to
\begin{align*}
\frac 12&\phi\left( \begin{tikzpicture}[scale=0.6,tcenter]
            \draw[ultra thick,black,l]
                (3,0.5) -- (0.5,0.5) -- (0,0);
            \draw[ultra thick,black,l,red pseudo]
                (0,0) -- (0,1);
            \draw[ultra thick,black,l]
                (0,0) -- (1,0)
                (0,1) -- (1,1);
            \draw[black,c=4pt]
                (0,0) node[fill] {}
                (0,1) node[fill] {}
                (0.5,0.5) node[fill] {};
            \draw[black,circle]
                (3,0.5) node[nb] {$w$};
             \draw[black,ellipse]
                (1.5,1) node[nbe] {$x\triangle y$}
				(1.5,0) node[nbe] {$y\triangle z$};
		\end{tikzpicture}\right)
-\frac 12\phi\left( \begin{tikzpicture}[scale=0.6,tcenter]
            \draw[ultra thick,black,l]
                (3,0.5) -- (0.5,0.5) -- (0,1)
                (0.5,0.5) -- (0,0);
            \draw[ultra thick,black,l,red pseudo]
                (0,0) -- (0,1);
            \draw[ultra thick,black,l]
                (0,0) -- (1,0)
                (0,1) -- (1,1);
            \draw[black,c=4pt]
                (0,0) node[fill] {}
                (0,1) node[fill] {}
                (0.5,0.5) node[fill] {};
             \draw[black,circle]
                (3,0.5) node[nb] {$w$};
             \draw[black,ellipse]
                (1.5,1) node[nbe] {$x\triangle y$}
				(1.5,0) node[nbe] {$y\triangle z$};
		\end{tikzpicture}\right)
-\frac 12 \phi\left( \begin{tikzpicture}[scale=0.6,tcenter]
            \draw[ultra thick,black,l]
                (3,0.5) -- (0.5,0.5) -- (0,1);
            \draw[ultra thick,black,l,red pseudo]
                (0,0) -- (0,1);
            \draw[ultra thick,black,l]
                (0,0) -- (1,0)
                (0,1) -- (1,1);
            \draw[black,c=4pt]
                (0,0) node[fill] {}
                (0,1) node[fill] {}
                (0.5,0.5) node[fill] {};
            \draw[black,circle]
                (3,0.5) node[nb] {$w$};
             \draw[black,ellipse]
                (1.5,1) node[nbe] {$x\triangle y$}
				(1.5,0) node[nbe] {$y\triangle z$};
		\end{tikzpicture}\right)
+ \frac 12	\phi\left( \begin{tikzpicture}[scale=0.6,tcenter]
            \draw[ultra thick,black,l]
                (3,0.5) -- (0.5,0.5);
            \draw[ultra thick,black,l,red pseudo]
                (0,0) -- (0,1);
            \draw[ultra thick,black,l]
                (0,0) -- (1,0)
                (0,1) -- (1,1);
            \draw[black,c=4pt]
                (0,0) node[fill] {}
                (0,1) node[fill] {}
                (0.5,0.5) node[fill] {};
            \draw[black,circle]
                (3,0.5) node[nb] {$w$};
             \draw[black,ellipse]
                (1.5,1) node[nbe] {$x\triangle y$}
				(1.5,0) node[nbe] {$y\triangle z$};
		\end{tikzpicture}\right)\\
=& \frac 18  \phi\left( \begin{tikzpicture}[scale=0.6,tcenter]
            \draw[ultra thick,black,l,red pseudo]
                (0,0) -- (0,1);
            \draw[ultra thick,black,l]
                (0,0) -- (1,0)
                (0,1) -- (1,1);
            \draw[black,c=4pt]
                (0,0) node[fill] {}
                (0,1) node[fill] {};
            \draw[black,ellipse]
                (1.5,1) node[nbe] {$x\triangle y$}
				(1.7,0) node[nbe = 33pt] {$y\triangle z\triangle w$};
		\end{tikzpicture}\right)
- \frac 18  \phi\left( \begin{tikzpicture}[scale=0.6,tcenter]
            \draw[ultra thick,black,l]
                (0,0) -- (1,0)
                (0,1) -- (1,1);
            \draw[black,c=4pt]
                (0,0) node[fill] {}
                (0,1) node[fill] {};
            \draw[black,ellipse]
                (1.5,1) node[nbe] {$x\triangle y$}
				(1.7,0) node[nbe = 33pt] {$y\triangle z\triangle w$};
		\end{tikzpicture}\right) 
+ \frac 18 \phi\left( \begin{tikzpicture}[scale=0.6,tcenter]
            \draw[ultra thick,black,l]
                (0,0) -- (1,0)
                (0,1) -- (1,1)
                (0,0) -- (0,1);
            \draw[black,c=4pt]
                (0,0) node[fill] {}
                (0,1) node[fill] {};
            \draw[black,ellipse]
                (1.5,1) node[nbe] {$x\triangle y$}
				(1.7,0) node[nbe = 33pt] {$y\triangle z \triangle w$};
		\end{tikzpicture}\right)
=\frac  14 		 \phi\left( \begin{tikzpicture}[scale=0.6,tcenter]
            \draw[ultra thick,black,l,red pseudo]
                (0,0) -- (0,1);
            \draw[ultra thick,black,l]
                (0,0) -- (1,0)
                (0,1) -- (1,1);
            \draw[black,c=4pt]
                (0,0) node[fill] {}
                (0,1) node[fill] {};
            \draw[black,ellipse]
                (1.5,1) node[nbe] {$x\triangle y$}
				(1.7,0) node[nbe = 33pt] {$y\triangle z\triangle w$};
		\end{tikzpicture}\right),
\end{align*}
Here, in order to obtain the second line from the first, we have  applied the deletion-contraction  relation~(\ref{rel:del-cont}) to the first and the fourth, and to the second and the third summands. The right-hand part equals:
\begin{align*}
- &\phi\left( \begin{tikzpicture}[scale=0.6,tcenter]
            \draw[ultra thick,black,l]
                (3,0.5) -- (0.5,0.5) -- (0,1)
                (0.5,0.5) -- (0,0);
            \draw[ultra thick,black,l,red pseudo]
                (0,0) -- (0,1);
            \draw[ultra thick,black,l]
                (0,0) -- (1,0)
                (0,1) -- (1,1);
            \draw[black,c=4pt]
                (0,0) node[fill] {}
                (0,1) node[fill] {}
                (0.5,0.5) node[fill] {};
             \draw[black,circle]
                (3,0.5) node[nb] {$w$};
             \draw[black,ellipse]
                (1.5,1) node[nbe] {$x\triangle y$}
				(1.5,0) node[nbe] {$y\triangle z$};
		\end{tikzpicture}\right)
+ \phi\left( \begin{tikzpicture}[scale=0.6,tcenter]
            \draw[ultra thick,black,l]
                (3,0.5) -- (0.5,0.5);
            \draw[ultra thick,black,l,red pseudo]
                (0,0) -- (0,1);
            \draw[ultra thick,black,l]
                (0,0) -- (1,0)
                (0,1) -- (1,1);
            \draw[black,c=4pt]
                (0,0) node[fill] {}
                (0,1) node[fill] {}
                (0.5,0.5) node[fill] {};
            \draw[black,circle]
                (3,0.5) node[nb] {$w$};
             \draw[black,ellipse]
                (1.5,1) node[nbe] {$x\triangle y$}
				(1.5,0) node[nbe] {$y\triangle z$};
		\end{tikzpicture}\right)
-\frac 14 \phi\left( \begin{tikzpicture}[scale=0.6,tcenter]
            \draw[ultra thick,black,l,red pseudo]
                (0,0) -- (0,1);
            \draw[ultra thick,black,l]
                (0,0) -- (1,0)
                (0,1) -- (1,1);
            \draw[black,c=4pt]
                (0,0) node[fill] {}
                (0,1) node[fill] {};
            \draw[black,ellipse]
                (1.7,1) node[nbe = 33 pt] {$x\triangle y\triangle w$}
				(1.5,0) node[nbe] {$y\triangle z$};
		\end{tikzpicture}\right) \\
= - &\phi\left( \begin{tikzpicture}[scale=0.6,tcenter]
            \draw[ultra thick,black,l]
                (3,0.5) -- (0.5,0.5) -- (0,1)
                (0.5,0.5) -- (0,0);
            \draw[ultra thick,black,l,red pseudo]
                (0,0) -- (0,1);
            \draw[ultra thick,black,l]
                (0,0) -- (1,0)
                (0,1) -- (1,1);
            \draw[black,c=4pt]
                (0,0) node[fill] {}
                (0,1) node[fill] {}
                (0.5,0.5) node[fill] {};
             \draw[black,circle]
                (3,0.5) node[nb] {$w$};
             \draw[black,ellipse]
                (1.5,1) node[nbe] {$x\triangle y$}
				(1.5,0) node[nbe] {$y\triangle z$};
		\end{tikzpicture}\right)
- \phi\left( \begin{tikzpicture}[scale=0.6,tcenter]
            \draw[ultra thick,black,l]
                (3,0.5) -- (0.5,0.5) -- (0,1);
            \draw[ultra thick,black,l,red pseudo]
                (0,0) -- (0,1);
            \draw[ultra thick,black,l]
                (0,0) -- (1,0)
                (0,1) -- (1,1);
            \draw[black,c=4pt]
                (0,0) node[fill] {}
                (0,1) node[fill] {}
                (0.5,0.5) node[fill] {};
            \draw[black,circle]
                (3,0.5) node[nb] {$w$};
             \draw[black,ellipse]
                (1.5,1) node[nbe] {$x\triangle y$}
				(1.5,0) node[nbe] {$y\triangle z$};
		\end{tikzpicture}\right)
= \frac 14  \phi\left( \begin{tikzpicture}[scale=0.6,tcenter]
            \draw[ultra thick,black,l,red pseudo]
                (0,0) -- (0,1);
            \draw[ultra thick,black,l]
                (0,0) -- (1,0)
                (0,1) -- (1,1);
            \draw[black,c=4pt]
                (0,0) node[fill] {}
                (0,1) node[fill] {};
            \draw[black,ellipse]
                (1.5,1) node[nbe] {$x\triangle y$}
				(1.7,0) node[nbe=33pt] {$y\triangle z\triangle w$};
		\end{tikzpicture}\right),
\end{align*}
where, again, relation~(\ref{rel:del-cont}) is used. The contributions are indeed equal. The symmetric case when $s$ is connected to the topmost and to the rightmost vertices of the triangle is treated similarly.

\item If $s$ is connected to all three vertices of the triangle, the weighted sum of the proper 3-colorings on the left-hand side is 0. Expressing the graph with dashed edges as a linear combination of graphs with solid edges, we have:
\begin{align*}
      0 =   \phi\left( \begin{tikzpicture}[scale=0.6,tcenter]
            \draw[ultra thick,black,l]
                (0,0.3) -- (1,0.8);
            \draw[ultra thick,black,l,red pseudo]
                (-0.5,0) -- (0.5,0) -- (0,0.8) -- (-0.5,0) -- (0,0.3) -- (0,0.8)
                (0,0.3) -- (0.5,0);
            \draw[ultra thick,black,l]
                (-0.5,0) -- (-1,-0.5)
                (0.5,0) -- (1,-0.5)
                (0,0.8) -- (0,1.6);
            \draw[black,c=4pt]
                (-0.5,0) node[fill] {}
                (0,0.8) node[fill] {}
                (0.5,0) node[fill] {}
                (0,0.3) node[fill] {};
            \draw[black,circle]
                (-1,-0.5) node[nb] {$x$}
                (0,1.6) node[nb] {$y$}
                (1,-0.5) node[nb] {$z$}
                (1,0.8) node[nb] {$w$};
        \end{tikzpicture}\right) = 
        \phi\left( \begin{tikzpicture}[scale=0.6,tcenter]
            \draw[ultra thick,black,l]
                (0,0.3) -- (1,0.8);
            \draw[ultra thick,black,l,red pseudo]
                (-0.5,0) -- (0.5,0) -- (0,0.8) -- (-0.5,0);
            \draw[ultra thick,black,l]
                (-0.5,0) -- (-1,-0.5)
                (0.5,0) -- (1,-0.5)
                (0,0.8) -- (0,1.6)
                (0,0.8) -- (0,0.3) -- (-0.5,0)
                (0,0.3) -- (0.5,0);
            \draw[black,c=4pt]
                (-0.5,0) node[fill] {}
                (0,0.8) node[fill] {}
                (0.5,0) node[fill] {}
                (0,0.3) node[fill] {};
            \draw[black,circle]
                (-1,-0.5) node[nb] {$x$}
                (0,1.6) node[nb] {$y$}
                (1,-0.5) node[nb] {$z$}
                (1,0.8) node[nb] {$w$};
        \end{tikzpicture}\right) -
         \phi\left( \begin{tikzpicture}[scale=0.6,tcenter]
           \draw[ultra thick,black,l]
                (0,0.3) -- (1,0.8);
            \draw[ultra thick,black,l,red pseudo]
                (-0.5,0) -- (0.5,0) -- (0,0.8) -- (-0.5,0);
            \draw[ultra thick,black,l]
                (-0.5,0) -- (-1,-0.5)
                (0.5,0) -- (1,-0.5)
                (0,0.8) -- (0,1.6)
                (0,0.8) -- (0,0.3) -- (-0.5,0);
            \draw[black,c=4pt]
                (-0.5,0) node[fill] {}
                (0,0.8) node[fill] {}
                (0.5,0) node[fill] {}
                (0,0.3) node[fill] {};
            \draw[black,circle]
                (-1,-0.5) node[nb] {$x$}
                (0,1.6) node[nb] {$y$}
                (1,-0.5) node[nb] {$z$}
                (1,0.8) node[nb] {$w$};
        \end{tikzpicture}\right)-
         \phi\left( \begin{tikzpicture}[scale=0.6,tcenter]
            \draw[ultra thick,black,l]
                (0,0.3) -- (1,0.8);
            \draw[ultra thick,black,l,red pseudo]
                (-0.5,0) -- (0.5,0) -- (0,0.8) -- (-0.5,0);
            \draw[ultra thick,black,l]
                (-0.5,0) -- (-1,-0.5)
                (0.5,0) -- (1,-0.5)
                (0,0.8) -- (0,1.6)
                (0,0.8) -- (0,0.3)
                (0,0.3) -- (0.5,0);
            \draw[black,c=4pt]
                (-0.5,0) node[fill] {}
                (0,0.8) node[fill] {}
                (0.5,0) node[fill] {}
                (0,0.3) node[fill] {};
            \draw[black,circle]
                (-1,-0.5) node[nb] {$x$}
                (0,1.6) node[nb] {$y$}
                (1,-0.5) node[nb] {$z$}
                (1,0.8) node[nb] {$w$};
        \end{tikzpicture}\right)-
         \phi\left( \begin{tikzpicture}[scale=0.6,tcenter]
           \draw[ultra thick,black,l]
                (0,0.3) -- (1,0.8);
            \draw[ultra thick,black,l,red pseudo]
                (-0.5,0) -- (0.5,0) -- (0,0.8) -- (-0.5,0);
            \draw[ultra thick,black,l]
                (-0.5,0) -- (-1,-0.5)
                (0.5,0) -- (1,-0.5)
                (0,0.8) -- (0,1.6)
                (0,0.3) -- (-0.5,0)
                (0,0.3) -- (0.5,0);
            \draw[black,c=4pt]
                (-0.5,0) node[fill] {}
                (0,0.8) node[fill] {}
                (0.5,0) node[fill] {}
                (0,0.3) node[fill] {};
            \draw[black,circle]
                (-1,-0.5) node[nb] {$x$}
                (0,1.6) node[nb] {$y$}
                (1,-0.5) node[nb] {$z$}
                (1,0.8) node[nb] {$w$};
        \end{tikzpicture}\right)\\
        + \phi\left( \begin{tikzpicture}[scale=0.6,tcenter]
           \draw[ultra thick,black,l]
                (0,0.3) -- (1,0.8);
            \draw[ultra thick,black,l,red pseudo]
                (-0.5,0) -- (0.5,0) -- (0,0.8) -- (-0.5,0);
            \draw[ultra thick,black,l]
                (-0.5,0) -- (-1,-0.5)
                (0.5,0) -- (1,-0.5)
                (0,0.8) -- (0,1.6)
                (0,0.8) -- (0,0.3);
            \draw[black,c=4pt]
                (-0.5,0) node[fill] {}
                (0,0.8) node[fill] {}
                (0.5,0) node[fill] {}
                (0,0.3) node[fill] {};
            \draw[black,circle]
                (-1,-0.5) node[nb] {$x$}
                (0,1.6) node[nb] {$y$}
                (1,-0.5) node[nb] {$z$}
                (1,0.8) node[nb] {$w$};
        \end{tikzpicture}\right)
         + \phi\left( \begin{tikzpicture}[scale=0.6,tcenter]
            \draw[ultra thick,black,l]
                (0,0.3) -- (1,0.8);
            \draw[ultra thick,black,l,red pseudo]
                (-0.5,0) -- (0.5,0) -- (0,0.8) -- (-0.5,0);
            \draw[ultra thick,black,l]
                (-0.5,0) -- (-1,-0.5)
                (0.5,0) -- (1,-0.5)
                (0,0.8) -- (0,1.6)
                (0.5,0) -- (0,0.3);
            \draw[black,c=4pt]
                (-0.5,0) node[fill] {}
                (0,0.8) node[fill] {}
                (0.5,0) node[fill] {}
                (0,0.3) node[fill] {};
            \draw[black,circle]
                (-1,-0.5) node[nb] {$x$}
                (0,1.6) node[nb] {$y$}
                (1,-0.5) node[nb] {$z$}
                (1,0.8) node[nb] {$w$};
        \end{tikzpicture}\right)
         + \phi\left( \begin{tikzpicture}[scale=0.6,tcenter]
            \draw[ultra thick,black,l]
                (0,0.3) -- (1,0.8);
            \draw[ultra thick,black,l,red pseudo]
                (-0.5,0) -- (0.5,0) -- (0,0.8) -- (-0.5,0);
            \draw[ultra thick,black,l]
                (-0.5,0) -- (-1,-0.5)
                (0.5,0) -- (1,-0.5)
                (0,0.8) -- (0,1.6)
                (-0.5,0) -- (0,0.3);
            \draw[black,c=4pt]
                (-0.5,0) node[fill] {}
                (0,0.8) node[fill] {}
                (0.5,0) node[fill] {}
                (0,0.3) node[fill] {};
            \draw[black,circle]
                (-1,-0.5) node[nb] {$x$}
                (0,1.6) node[nb] {$y$}
                (1,-0.5) node[nb] {$z$}
                (1,0.8) node[nb] {$w$};
        \end{tikzpicture}\right)
         - \phi\left( \begin{tikzpicture}[scale=0.6,tcenter]
            \draw[ultra thick,black,l]
                (0,0.3) -- (1,0.8);
            \draw[ultra thick,black,l,red pseudo]
                (-0.5,0) -- (0.5,0) -- (0,0.8) -- (-0.5,0);
            \draw[ultra thick,black,l]
                (-0.5,0) -- (-1,-0.5)
                (0.5,0) -- (1,-0.5)
                (0,0.8) -- (0,1.6);
            \draw[black,c=4pt]
                (-0.5,0) node[fill] {}
                (0,0.8) node[fill] {}
                (0.5,0) node[fill] {}
                (0,0.3) node[fill] {};
            \draw[black,circle]
                (-1,-0.5) node[nb] {$x$}
                (0,1.6) node[nb] {$y$}
                (1,-0.5) node[nb] {$z$}
                (1,0.8) node[nb] {$w$};
        \end{tikzpicture}\right).
\end{align*}
The corresponding contribution to the right-hand part of equality~(\ref{rel:triangle}), up to the overall factor of $\frac 12$, should be equal to
\begin{align*}
\phi\left( \begin{tikzpicture}[scale=0.6,tcenter]
            \draw[ultra thick,black,l]
                (0.5,0.5) -- (3,0.5);
            \draw[ultra thick,black,l,red pseudo]
                (0,0) -- (0,1);
            \draw[ultra thick,black,l]
                (0,0) -- (1,0)
                (0,1) -- (1,1);
            \draw[black,c=4pt]
                (0,0) node[fill] {}
                (0,1) node[fill] {}
                (0.5,0.5) node[fill] {};
            \draw[black,circle]
                (3,0.5) node[nb] {$w$};
             \draw[black,ellipse]
                (1.5,1) node[nbe] {$x\triangle y$}
				(1.5,0) node[nbe] {$y\triangle z$};
		\end{tikzpicture}\right)
  -\phi\left( \begin{tikzpicture}[scale=0.6,tcenter]
            \draw[ultra thick,black,l]
                (3,0.5) -- (0.5,0.5) -- (0,0);
            \draw[ultra thick,black,l,red pseudo]
                (0,0) -- (0,1);
            \draw[ultra thick,black,l]
                (0,0) -- (1,0)
                (0,1) -- (1,1);
            \draw[black,c=4pt]
                (0,0) node[fill] {}
                (0,1) node[fill] {}
                (0.5,0.5) node[fill] {};
            \draw[black,circle]
                (3,0.5) node[nb] {$w$};
             \draw[black,ellipse]
                (1.5,1) node[nbe] {$x\triangle y$}
				(1.5,0) node[nbe] {$y\triangle z$};
		\end{tikzpicture}\right)
 - \phi\left( \begin{tikzpicture}[scale=0.6,tcenter]
            \draw[ultra thick,black,l]
                (3,0.5) -- (0.5,0.5) -- (0,1);
            \draw[ultra thick,black,l,red pseudo]
                (0,0) -- (0,1);
            \draw[ultra thick,black,l]
                (0,0) -- (1,0)
                (0,1) -- (1,1);
            \draw[black,c=4pt]
                (0,0) node[fill] {}
                (0,1) node[fill] {}
                (0.5,0.5) node[fill] {};
            \draw[black,circle]
                (3,0.5) node[nb] {$w$};
             \draw[black,ellipse]
                (1.5,1) node[nbe] {$x\triangle y$}
				(1.5,0) node[nbe] {$y\triangle z$};
		\end{tikzpicture}\right)
   - \phi\left( \begin{tikzpicture}[scale=0.6,tcenter]
            \draw[ultra thick,black,l]
                (3,0.5) -- (0.5,0.5) -- (0,1)
                (0.5,0.5) -- (0,0);
            \draw[ultra thick,black,l,red pseudo]
                (0,0) -- (0,1);
            \draw[ultra thick,black,l]
                (0,0) -- (1,0)
                (0,1) -- (1,1);
            \draw[black,c=4pt]
                (0,0) node[fill] {}
                (0,1) node[fill] {}
                (0.5,0.5) node[fill] {};
             \draw[black,circle]
                (3,0.5) node[nb] {$w$};
             \draw[black,ellipse]
                (1.5,1) node[nbe] {$x\triangle y$}
				(1.5,0) node[nbe] {$y\triangle z$};
		\end{tikzpicture}\right)\\
  +  \phi\left( \begin{tikzpicture}[scale=0.6,tcenter]
            \draw[ultra thick,black,l]
                (3,0.5) -- (0.5,0.5) -- (0,1)
                (0.5,0.5) -- (0,0);
            \draw[ultra thick,black,l,red pseudo]
                (0,0) -- (0,1);
            \draw[ultra thick,black,l]
                (0,0) -- (1,0)
                (0,1) -- (1,1);
            \draw[black,c=4pt]
                (0,0) node[fill] {}
                (0,1) node[fill] {}
                (0.5,0.5) node[fill] {};
              \draw[black,circle]
                (3,0.5) node[nb] {$w$};
             \draw[black,ellipse]
                (1.5,1) node[nbe] {$x\triangle y$}
				(1.5,0) node[nbe] {$y\triangle z$};
		\end{tikzpicture}\right)
   + \phi\left( \begin{tikzpicture}[scale=0.6,tcenter]
            \draw[ultra thick,black,l]
                (3,0.5) -- (0.5,0.5) -- (0,0);
            \draw[ultra thick,black,l,red pseudo]
                (0,0) -- (0,1);
            \draw[ultra thick,black,l]
                (0,0) -- (1,0)
                (0,1) -- (1,1);
            \draw[black,c=4pt]
                (0,0) node[fill] {}
                (0,1) node[fill] {}
                (0.5,0.5) node[fill] {};
               \draw[black,circle]
                (3,0.5) node[nb] {$w$};
             \draw[black,ellipse]
                (1.5,1) node[nbe] {$x\triangle y$}
				(1.5,0) node[nbe] {$y\triangle z$};
		\end{tikzpicture}\right)
   + \phi\left( \begin{tikzpicture}[scale=0.6,tcenter]
            \draw[ultra thick,black,l]
                (3,0.5) -- (0.5,0.5) -- (0,1);
            \draw[ultra thick,black,l,red pseudo]
                (0,0) -- (0,1);
            \draw[ultra thick,black,l]
                (0,0) -- (1,0)
                (0,1) -- (1,1);
            \draw[black,c=4pt]
                (0,0) node[fill] {}
                (0,1) node[fill] {}
                (0.5,0.5) node[fill] {};
            \draw[black,circle]
                (3,0.5) node[nb] {$w$};
             \draw[black,ellipse]
                (1.5,1) node[nbe] {$x\triangle y$}
				(1.5,0) node[nbe] {$y\triangle z$};
		\end{tikzpicture}\right)
   - \phi\left( \begin{tikzpicture}[scale=0.6,tcenter]
            \draw[ultra thick,black,l]
                (3,0.5) -- (0.5,0.5);
            \draw[ultra thick,black,l,red pseudo]
                (0,0) -- (0,1);
            \draw[ultra thick,black,l]
                (0,0) -- (1,0)
                (0,1) -- (1,1);
            \draw[black,c=4pt]
                (0,0) node[fill] {}
                (0,1) node[fill] {}
                (0.5,0.5) node[fill] {};
            \draw[black,circle]
                (3,0.5) node[nb] {$w$};
             \draw[black,ellipse]
                (1.5,1) node[nbe] {$x\triangle y$}
				(1.5,0) node[nbe] {$y\triangle z$};
		\end{tikzpicture}\right)
    \end{align*}
and vanish identically as well.
\end{itemize}
The claim follows.
\end{proof}

Notice that the graph Chmutov-Varchenko relations~(\ref{defn:ChV}) imply the triangle identity~(\ref{rel:triangle}), as it is the difference of the second and the first graph 6T-relations. It turns out that the triangle identity implies the 4T-relations~(\ref{rel:gr4T}):

\begin{thm}\label{thm:4inv}
    Function $\phi$ is a 4-invariant.
\end{thm}
\begin{proof}
    Indeed:
\begin{align*}
\phi\left( \begin{tikzpicture}[scale=0.6,tcenter]
            \draw[ultra thick,black,l,red pseudo]
                (0,0) -- (0,1);
            \draw[ultra thick,black,l]
                (0,0) -- (1,0)
                (0,1) -- (1,1);
            \draw[black,c=4pt]
                (0,0) node[fill] {}
                (0,1) node[fill] {};
            \draw[black,circle]
                (1.0,1) node[nb] {$a$}
				(1.0,0) node[nb] {$b$};
		\end{tikzpicture}\right)
= 2 \phi\left( \begin{tikzpicture}[scale=0.6,tcenter]
            \draw[ultra thick,black,l,red pseudo]
                (-0.5,0) -- (0.5,0) -- (0,0.8) -- (-0.5,0);
            \draw[ultra thick,black,l]
                (-0.5,0) -- (-1,-0.5)
                (0.5,0) -- (1,-0.5);
            \draw[black,c=4pt]
                (-0.5,0) node[fill] {}
                (0,0.8) node[fill] {}
                (0.5,0) node[fill] {};
            \draw[black,circle]
                (-1,-0.5) node[nb] {$a$}
                (1,-0.5) node[nb] {$b$};
        \end{tikzpicture}\right)
= 2 \phi\left( \begin{tikzpicture}[scale=0.6,tcenter]
            \draw[ultra thick,black,l,red pseudo]
                (-0.5,0) -- (0.5,0) -- (0,0.8) -- (-0.5,0);
            \draw[ultra thick,black,l]
                (-0.5,0) -- (-1,-0.5)
                (0,0.8) -- (0,1.6);
            \draw[black,c=4pt]
                (-0.5,0) node[fill] {}
                (0,0.8) node[fill] {}
                (0.5,0) node[fill] {};
            \draw[black,circle]
                (0,1.6) node[nb] {$a$}
                (-1,-0.5) node[nb] {$b$};
        \end{tikzpicture}\right)
=\phi\left( \begin{tikzpicture}[scale=0.6,tcenter]
            \draw[ultra thick,black,l,red pseudo]
                (0,0) -- (0,1);
            \draw[ultra thick,black,l]
                (0,0) -- (1.5,0)
                (0,1) -- (1,1);
            \draw[black,c=4pt]
                (0,0) node[fill] {}
                (0,1) node[fill] {};
            \draw[black,circle]
				(1.5,0) node[nb] {$a$};
			\draw[black,ellipse]
                (1.5,1) node[nbe] {$a\triangle b$};
		\end{tikzpicture}\right).
\end{align*}
\end{proof}

\subsection{The graph 6T-relations for \texorpdfstring{$\phi$}{φ}}

\begin{prop}\label{prop:6T}
Function $\phi$ satisfies the graph 6T-relations of definition~\ref{defn:ChV}.
\end{prop}
\begin{proof}
As it was proven in the previous subsection, function $\phi$ satisfies the triangle identity~(\ref{rel:triangle}), and hence the difference of  6T-relation holds:
\begin{align*}
 \phi \left( \begin{tikzpicture}[scale=0.4,tcenter]
                    \draw[ultra thick, black, l]
                        (0,1) -- (1.5,1)
                        (0,0) -- (1.5,0)
                        (0,-1) -- (1.5,-1)
                        (0,1) edge[out=-120,in=120] (0,-1);
                    \draw[ultra thick, red pseudo, l]
                        (0,1) -- (0,-1);
                    \draw[black, c=4pt]
                        (0,1) node[fill] {}
                        (0,0) node[fill] {}
                        (0,-1) node[fill] {};
                    \draw[black,circle]
                        (1.5,1) node[nb] {$x$}
                        (1.5,0) node[nb] {$y$}
                        (1.5,-1) node[nb] {$z$};
                \end{tikzpicture}\right) 
-\phi\left( \begin{tikzpicture}[scale=0.4,tcenter]
                    \draw[ultra thick, red pseudo, l]
                        (0,1) -- (0,-1);
                    \draw[ultra thick, black, l]
                        (0,1) -- (1.5,1)
                        (0,0) -- (1.5,0)
                        (0,-1) -- (1.5,-1);
                    \draw[black, c=4pt]
                        (0,1) node[fill] {}
                        (0,0) node[fill] {}
                        (0,-1) node[fill] {};
                    \draw[black,circle]
                        (1.5,1) node[nb] {$x$}
                        (1.5,0) node[nb] {$y$}
                        (1.5,-1) node[nb] {$z$};
                \end{tikzpicture}\right) 
                =   \phi\left( \begin{tikzpicture}[scale=0.6,tcenter]
            \draw[ultra thick,black,l,red pseudo]
                (-0.5,0) -- (0.5,0) -- (0,0.8) -- (-0.5,0);
            \draw[ultra thick,black,l]
                (-0.5,0) -- (-1,-0.5)
                (0.5,0) -- (1,-0.5)
                (0,0.8) -- (0,1.6);
            \draw[black,c=4pt]
                (-0.5,0) node[fill] {}
                (0,0.8) node[fill] {}
                (0.5,0) node[fill] {};
            \draw[black,circle]
                (-1,-0.5) node[nb] {$x$}
                (0,1.6) node[nb] {$y$}
                (1,-0.5) node[nb] {$z$};
        \end{tikzpicture}\right)
        =\frac{1}{2}\phi\left( \begin{tikzpicture}[scale=0.6,tcenter]
            \draw[ultra thick,black,l,red pseudo]
                (0,0) -- (0,1);
            \draw[ultra thick,black,l]
                (0,0) -- (1,0)
                (0,1) -- (1,1);
            \draw[black,c=4pt]
                (0,0) node[fill] {}
                (0,1) node[fill] {};
            \draw[black,ellipse]
                (1.5,1) node[nbe] {$x\triangle y$}
				(1.5,0) node[nbe] {$y\triangle z$};
		\end{tikzpicture}\right).
\end{align*}

So it remains to prove that the sum of the graph 6T-relations holds for $\phi$. On the left-hand side, due to relation~(\ref{rel:del-cont-var}) and the graph 4T-relations of definition~\ref{rel:4T}, we have:
\begin{align*}
 \phi \left( \begin{tikzpicture}[scale=0.4,tcenter]
                    \draw[ultra thick, black, l]
                        (0,1) -- (1.5,1)
                        (0,0) -- (1.5,0)
                        (0,-1) -- (1.5,-1)
                        (0,1) edge[out=-120,in=120] (0,-1);
                    \draw[ultra thick, red pseudo, l]
                        (0,1) -- (0,-1);
                    \draw[black, c=4pt]
                        (0,1) node[fill] {}
                        (0,0) node[fill] {}
                        (0,-1) node[fill] {};
                    \draw[black,circle]
                        (1.5,1) node[nb] {$x$}
                        (1.5,0) node[nb] {$y$}
                        (1.5,-1) node[nb] {$z$};
                \end{tikzpicture}\right) 
+\phi\left( \begin{tikzpicture}[scale=0.4,tcenter]
                    \draw[ultra thick, red pseudo, l]
                        (0,1) -- (0,-1);
                    \draw[ultra thick, black, l]
                        (0,1) -- (1.5,1)
                        (0,0) -- (1.5,0)
                        (0,-1) -- (1.5,-1);
                    \draw[black, c=4pt]
                        (0,1) node[fill] {}
                        (0,0) node[fill] {}
                        (0,-1) node[fill] {};
                    \draw[black,circle]
                        (1.5,1) node[nb] {$x$}
                        (1.5,0) node[nb] {$y$}
                        (1.5,-1) node[nb] {$z$};
                \end{tikzpicture}\right) 
= -\frac{1}{2} \phi\left( \begin{tikzpicture}[scale=0.6,tcenter]
            \draw[ultra thick,red pseudo,l]
                (0,0) -- (0,1);
            \draw[ultra thick,black,l]
                (0,0) -- (1.5,0)
                (0,1) -- (1.5,1);
            \draw[black,c=4pt]
                (0,0) node[fill] {}
                (0,1) node[fill] {};
            \draw[black,circle]
                (1.5,1) node[nb] {$y$};
            \draw[black,ellipse]
                (1.5,0) node[nbe] {$x\triangle z$};
        \end{tikzpicture}\right)
=-\frac{1}{2} \phi\left( \begin{tikzpicture}[scale=0.6,tcenter]
            \draw[ultra thick,red pseudo,l]
                (0,0) -- (0,1);
            \draw[ultra thick,black,l]
                (0,0) -- (1.5,0)
                (0,1) -- (1.5,1);
            \draw[black,c=4pt]
                (0,0) node[fill] {}
                (0,1) node[fill] {};
            \draw[black,circle]
                (1.5,1) node[nb] {$y$};
            \draw[black,ellipse]
                (1.7,0) node[nbe=33pt] {$x\triangle y \triangle z$};
        \end{tikzpicture}\right).
\end{align*}
While on the left-hand side of the sum we get:
\begin{align*}
\phi\left( \begin{tikzpicture}[scale=0.4,tcenter]
                    \draw[ultra thick, black, l]
                        (0,1) -- (3,1)
                        (3,-1) -- (0,-1);
                    \draw[black, c=4pt]
                        (0,1) node[fill] {}
                        (0,-1) node[fill] {};
                    \draw[black,ellipse]
                        (3,1) node[nb] {$y$}
                        (3,-1) node[nbe=33pt] {$x\triangle y\triangle z$};
                \end{tikzpicture}\right)
-\frac{1}{2} \phi \left( \begin{tikzpicture}[scale=0.4,tcenter]
                    \draw[ultra thick, black, l]
                        (0,1) -- (1.5,1)
                        (1.5,-1) -- (0,-1);
                    \draw[black, c=4pt]
                        (0,1) node[fill] {}
                        (0,-1) node[fill] {};
                    \draw[black,ellipse]
                        (2,1) node[nbe] {$x\triangle y$}
                        (2,-1) node[nbe] {$y\triangle z$};
                \end{tikzpicture}\right)
-\frac{1}{2} \phi \left( \begin{tikzpicture}[scale=0.4,tcenter]
                    \draw[ultra thick, black, l]
                        (1.5,-1) -- (0,-1) -- (0,1) -- (1.5,1);
                    \draw[black, c=4pt]
                        (0,1) node[fill] {}
                        (0,-1) node[fill] {};
                    \draw[black,ellipse]
                        (2,1) node[nbe] {$x\triangle y$}
                        (2,-1) node[nbe] {$y\triangle z$};
                \end{tikzpicture}\right).
\end{align*}
Due to the deletion-contraction relation~(\ref{rel:del-cont}), these expressions are equal.

\end{proof}

Combining together propositions~\ref{prop:mulnorm}, \ref{prop:ldel}, \ref{prop:6T}, and theorem~\ref{thm:4inv} we get:
\begin{thm}\label{thm: main}
The 4-invariant $\phi$ satisfies the graph Chmutov-Varchenko relations at $\frac 38$ and thus extends the value of $w_{\mathfrak{sl}(2)}$ at $c = \frac 38$ from the set of graphs realizable as intersection graphs of chord diagrams to the whole set $\mathbf G$.
\end{thm}

\begin{remark}
	The graph Chmutov-Varchenko relations at $p\ne 0,\frac 38$ do not produce a well-defined graph invariant.
\end{remark}

\begin{remark}\label{remark:sl2-2d}
	After the first version of this text was written, P. Zakorko brought to our attention that the specialization of the universal $\mathfrak{sl}(2)$-weight system to the irreducible 2-dimensional representation of $\mathfrak{sl}(2)$ satisfies a deletion-contraction relation analogous to the one in~(\ref{rel:del-cont}).
	Additionally, it assigns the same value to the chord diagram with a single chord as the 4-invariant $\phi$.
	That implies that $\phi$ can be viewed as an extension of the universal $\mathfrak{sl}(2)$-weight system specialized to the standard 2-dimensional representation.
	The value $\frac{3}{8}$, which $\phi$ assigns to the single-vertex graph, corresponds to the eigenvalue of the Casimir element in the 2-dimensional representation of $\sl(2)$.

	It is known (see, e.g.,~\cite{CDBook}) that the value of the $\sl(2)$-weight system on its 2-dimensional irreducible representation can be explicitly computed for any chord diagram.
	Moreover, this formula can be extended to graphs as follows:

	\begin{align*}
		\psi:
		G \mapsto
		\frac{1}{2^{2|V(G)|}} \sum_{U\subset V(G)} \left(-\frac{1}{2}\right)^{|V(G)|-|U|} 2^{\corank(A(G|_U))},
	\end{align*}
	where
	\begin{itemize}
		\item the sum is taken over all possible subsets $U$ of the vertex set $V(G)$;
		\item $G|_U$ is the subgraph of $G$ induced by vertex subset $U$;
		\item $A(G)$ is the adjacency matrix of $G$ with coefficients in the field with two elements;
		\item $\corank$ is the corank of the matrix.
	\end{itemize}

	However, the fact that this expression defines a 4-invariant (it follows, for instance, from the fact that it can be expressed as a convolution product of two 4-invariants) has never been explicitly stated in the literature to the best of our knowledge.

	We conjecture that $\phi = \psi$.
\end{remark}

\section{Further discussion}

\subsection{Alternative formulae for \texorpdfstring{$\phi$}{ φ }}
Here we discuss an alternative formula for $\phi$ that potentially allows for a faster computation of $\phi$.
\begin{prop}\label{prop:alt1}
	For a graph $G\in \mathbf{G}$ with the set of vertices $V(G)$ and the set of edges $E(G)$, we have
	\[
		\phi(G) = 2^{-3|V(G)|} \sum_{\substack{U \subset V(G)\\ G|_U \text{ is Eulerian}}} (-1)^{|E(U,V(G)\setminus U)|} 2^{|U|}.
	\]
	Here $G|_U$ is the full subgraph of $G$ induced by the subset of vertices $U$.
	A graph is Eulerian if all vertex degrees are even.
	The symbol $E(U, V(G) \setminus U)$ denotes the set of edges that connect vertices of $U$ to vertices of its complement.
\end{prop}
\begin{proof}
	\begin{align*}
		\phi(G)
		&= 2^{-3|V(G)|} \sum_{E'\subset E(G)} (-2)^{|E'|} \chi_3(G|_{E'}) = \\
		&= 2^{-3|V(G)|} \sum_{E'\subset E(G)} (-2)^{|E'|} \sum_{V(G)=V_1\sqcup V_2\sqcup V_3} \mathds{1}(V_1\sqcup V_2\sqcup V_3\text{ is a coloring of }G|_{E'}) = \\
		&= 2^{-3|V(G)|} \sum_{V(G)=V_1\sqcup V_2\sqcup V_3} \sum_{E'\subset E(G)} (-2)^{|E'|} \mathds{1}(V_1\sqcup V_2\sqcup V_3\text{ is a coloring of }G|_{E'}) = \\
		&= 2^{-3|V(G)|} \sum_{V(G)=V_1\sqcup V_2\sqcup V_3} \sum_{E'\subset E(V_1, V_2)\cup E(V_2, V_3)\cup E(V_3, V_1)} (-2)^{|E'|} = \\
		&= 2^{-3|V(G)|} \sum_{V(G)=V_1\sqcup V_2\sqcup V_3} \left(\sum_{E'\subset E(V_1, V_2)} (-2)^{|E'|}\right) \cdot \left(\sum_{E'\subset E(V_2, V_3)} (-2)^{|E'|}\right) \cdot \left(\sum_{E'\subset E(V_3, V_1)} (-2)^{|E'|}\right) = \\
		&= 2^{-3|V(G)|} \sum_{V(G)=V_1\sqcup V_2\sqcup V_3} (-1)^{|E(V_1, V_2)| + |E(V_2, V_3)| + |E(V_3, V_1)|} = \\
		&= 2^{-3|V(G)|} (-1)^{|E(G)|} \sum_{V(G)=V_1\sqcup V_2\sqcup V_3} (-1)^{|E(V_1, V_1)| + |E(V_2, V_2)| + |E(V_3, V_3)|} = \\
		&= 2^{-3|V(G)|} (-1)^{|E(G)|} \sum_{V(G)=V_1\sqcup V_2} (-1)^{|E(V_1, V_1)|} \sum_{V_2=V'_2\sqcup V''_2} (-1)^{|E(V'_2, V'_2)| + |E(V''_3, V''_3)|} = \\
		&= 2^{-3|V(G)|} (-1)^{|E(G)|} \sum_{V(G)=V_1\sqcup V_2} (-1)^{|E(V_1, V_1)| + |E(V_2, V_2)|} \sum_{V_2=V'_2\sqcup V''_2} (-1)^{|E(V'_2, V''_2)|} = \\
		&= 2^{-3|V(G)|} \sum_{V(G)=V_1\sqcup V_2} (-1)^{|E(V_1, V_2)|} \sum_{V_2=V'_2\sqcup V''_2} (-1)^{|E(V'_2, V''_2)|}.
	\end{align*}

	Here $\mathds{1}(\ldots)$ equals $1$ if its argument proposition is true and $0$ otherwise.

	To finish the proof, we will show that $\sum_{V(G)=V_1\sqcup V_2} (-1)^{|E(V_1, V_2)|} = 2^{|V(G)|} \mathds{1}(G\text{ is Eulerian})$.

	Let $G$ be Eulerian: each connected component admits an Eulerian cycle.
	For each $V(G)=V_1\sqcup V_2$, each component therefore contributes an even number of edges to $|E(V_1, V_2)|$, hence the sum simplifies to $\sum_{V(G)=V_1\sqcup V_2} 1 = 2^{|V(G)|}$.

	Let $G$ be non-Eulerian: there exists a vertex $v$ such that $\deg v$ is odd.
	For each $V(G)=V_1\sqcup V_2$, we consider $V(G)=V'_1\sqcup V'_2$, where $V'_i = V_i \triangle \setp{v}$.
	It is evident that $(-1)^{|E(V_1, V_2)|} + (-1)^{|E(V'_1,V'_2)|} = 0$.
	Since this correspondence is bijective, the sum evaluates to zero.
\end{proof}

Another expression comes from an explicit formula for the number of proper vertex 3-colorings.

\begin{prop}\label{def:phi-alt}
For a graph $G \in \mathbf G$ with the set of vertices $V(G)$ and the set of edges $E(G)$ we have
\[
\phi(G) = 2^{-3|V(G)|}(-1)^{|E(G)|} \sum_{E' \subset E(G)} (-2)^{|E'|} 3^{c(G|_{E'})}, 
\]
where $c(G|_{E'})$ is the number of the connected components of the graph obtained from $G$ by removing all the edges that do not belong to $E'$.
\end{prop}
\begin{proof}
Recall the alternative formula for the number of proper 3-colorings of the vertices of the graph that follows from the inclusion-exclusion argument (see, e.g.,~\cite{dong2005chromatic}):
\[
\chi_3(G) = \sum_{E' \subset E(G)} (-1)^{|E'|} 3^{c(G|_{E'})}
\]
and plug it into the definition of $\phi$:
\begin{align*}
\phi(G) = 2^{-3|V(G)|} \sum_{E' \subset E(G)} (-2)^{|E'|} \sum_{E'' \subset E'} (-1)^{E''} 3^{c(G|_{E''})}\\ =
2^{-3|V(G)|}\sum_{E'' \subset E(G)} 2^{|E''|} 3^{c(G|_{E''})}\sum_{E'' \subset E' \subset E(G)} (-2)^{|E'| - |E''|},
\end{align*}
where we change the order of summation: first sum over all possible choices of $E'' \subset E(G)$, and then over all the intermediate subsets $E'' \subset E' \subset E$. The last sum clearly equals $(-1)^{|E(G)| - |E''|}$. The assertion follows.
\end{proof}

\subsection{An extremal value of \texorpdfstring{$\phi$}{ φ } for the fixed number of vertices} 
%Notice, that the formula from proposition~\ref{def:phi-alt} implies the following algorithm for computation of $\phi$. Let $G\in \mathbf G_n$ and fix some order on the set of its vertices $v_1,\ldots,v_n$. For $k = 1,\ldots, n$ denote $G_k$ the graph induced by the first $k$ vertices. Let $P_k$ be the set of edges adjacent to the vertex $v_k$ in the graph $G_k$, and $E(G_k)$ -- the full set of edges of $G_k$. Then we have:
%\begin{align*}
%\phi(G_1) &= \frac 38,\\
%\phi(G_k) &= 2^{-3k}(-1)^{|P_k| + |E(G_{k-1})|}\left(3 \sum_{E' \subset E(G_{k-1})} (-2)^{|E'|}3^{c(G_{k-1}|_{E'})} + \sum_{\begin{smallmatrix}F \subset P_k\\ F\ne \emptyset \end{smallmatrix}} \sum_{E' \subset E(G_{k-1})} (-2)^{|E'|}3^{c(G_{k-1}|_{E'})}\right)\\
%&= \frac{(-1)^{|P_k|}\left(3 + (-1)^{|P_k|} - 1 \right)}8\phi(G_{k-1}).
%\end{align*} 
%In short, 
%\[\phi(G_k) = \begin{cases} \frac 38, &\mbox{ if $k = 1$ },\\
%\frac 38 \phi(G_{k-1}),&\mbox{ if $k>1$ and $|P_k|$ is even},\\
%-\frac 18 \phi(G_{k-1}),&\mbox{ if $k>1$ and $|P_k|$ is odd}.
%\end{cases}
%\] 
%In particular, this procedure allows to compute the values of $\phi$ on the complete graph $K_n$ on $n$ vertices:
%\[ \phi(K_n) = \frac {(-1)^{\lfloor \frac n2 \rfloor}3^{\lceil \frac n2 \rceil} }{8^n}.\]

We finish the paper with a proof of the following proposition motivated by an observation of M.~Kazarian and M.~Shapiro:
\begin{prop}
	Let $n \in \N$. For graph $G \in \mathbf{G}_n$, we have $\phi(G) \ne 0$ and
	\[
%		3\left( \frac{1}{8} \right)^{n}\le |\phi(G)| \le \left (\frac{3}{8} \right)^n.
		|\phi(G)| \le \left (\frac{3}{8} \right)^n.
	\]
	The bound is exact.
\end{prop}
\begin{proof}
	The claim follows easily from the formula from proposition~\ref{prop:alt1} with the sum over vertex subsets inducing Eulerian subgraphs. Namely
	\begin{align*}
|\phi(G)| \le 2^{-3|V(G)|}\sum_{U \subset V(G)} 2^{|U|} = \left( \frac 38 \right)^{|V(G)|}.
	\end{align*}
	
	The upper bound is achieved on discrete graphs.
\end{proof}

\bibliography{bib}

\end{document}